\documentclass[10pt]{amsart}

\usepackage{graphicx}
\usepackage{amsmath}
\usepackage{amssymb}
\usepackage{amsthm}
\usepackage{algorithm}
\usepackage{algorithmic}
\usepackage[top=26mm, bottom=26mm, left=26mm, right=26mm]{geometry}
\usepackage[dvips]{color}
\usepackage{bm}
\usepackage{marvosym}

\newcommand{\figref}[1]{{Figure~\ref{#1}}}
\newcommand{\tabref}[1]{{Table~\ref{#1}}}

\newcommand{\pd}[2]{\dfrac{\partial #1}{\partial #2}}

\begin{document}

\begin{center}
A fast topology optimisation for material- and geometry-independent cloaking devices
with the BEM and the ${\mathcal H}$-matrix method

\vspace{10pt}
Kenta Nakamoto$^1$, Hiroshi Isakari$^1$, Toru Takahashi$^1$, and Toshiro Matsumoto$^1$ \\
1. Nagoya University, Japan \\
\end{center}
\vspace{10pt}

\noindent
\textbf{Abstract:}
We show a design method of cloaking devices which
work for target objects with arbitrary shape and material
by a topology optimisation with an accurate and efficient sensitivity
analysis.
Most of past researches on topology optimisation of cloaking devices
intend to hide a circle-shaped perfect electric conductor. In this case, 
the cloaking effect is highly dependent on the shape and material 
of a target object. In this study, we consider to design a
cloaking device which work regardless of the property of target objects
by modifying the definition of the objective function.
Also, we developed an efficient and accurate sensitivity analysis with
the boundary element method and the ${\mathcal H}$-matrix method. We show that
the proposed method can successfully obtain desired cloaking devices
with low computational cost.

\vspace{10pt}

\newcommand{\bs}[1]{\boldsymbol{#1}}

\section{Introduction}
Establishment of the theory to control the electro-magnetic field with
metamaterials \cite{Pendry1999negative}
has encouraged developments of innovative optical devices such
as metamaterial-based waveguides \cite{Hrabar2005waveguide} and superlens
\cite{Fang2003superlens}, etc.
As one of the most attractive applications of metamaterials, cloaking
has extensively been researched so far.
The cloaking is a technology to make an object invisible by
reducing the scattering around the object.
There are two main ways to construct cloaking devices; light pass
controlling and scattering cancellation.
The former strategy was proposed by Pendry et al 
\cite{Pendry2006cloak},
in which they theoretically showed the possibility to construct a cloaking
device by periodically putting metamaterials whose permeability is negative.
By covering a target object with the cloaking device,

the light pass is controlled to go around the target object, 
which makes the object invisible.
Smith et al. experimentally created a cloaking device by allocating
metamaterials periodically and showed that their cloaking works
to hide a copper cylinder in TM polarised microwave
range\cite{Schurig2006}.
While their cloaking device can hide arbitrary-shaped 
target objects, it is made of anisotropic and inhomogeneous
materials which may cause difficulty in manufacturing.
Furthermore, their cloaking device has some critical problems from
an engineering standpoint. Their cloaking may cause an energy loss
and works only when the frequency of the incident field is in a narrow
target range.
As the latter type of cloaking, 
Al\'{u} and Engheta proposed to reduce scattering from
spherical and cylindrical objects by covering them with 
properly designed epsilon-near-zero (ENZ) material \cite{Alu2005} which
is a kind of metamaterials whose real part of macroscopic permittivity
is close to zero.
By the use of the ENZ material as a cover, we can suppress the dipolar
term in the Mie expansion of the scattered field from a plane wave.
Their cloaking device can be made of low loss isotropic and homogeneous
material and works in a wide bandwidth. 
Many researches on cloaking based on scattering cancellation with ENZ
material are reported \cite{Monti2011, Farhat2012}.
These cloaking devices, however, need to be designed carefully in accordance
with the shape of the target object. When the shape of the target object
is perturbated from the expected one, the cloaking may not work.
One can find that both methods have some pros and cons, and it is important
to develop a design method which is loss less, wide band, easily fabricated and
robust to change of shape of target object.

Recently, simulation-based design methods for engineering
devices have been applied to many varieties of fields.
Topology optimisation has extensively been researched as one of such
design methods.
Compared with the other methods such as
shape optimisation, topology optimisation is the most flexible design
method since we can determine not only peripheral shape but also
topology of devices.
The first pioneering research on an application of the topology
optimisation for engineering devices was conducted
by Bends\o e and Kikuchi in a design problem of elastic member to
maximise its stiffness \cite{Bendsoe1988}.
They proposed to express structure of devices as distribution of a
characteristic function and determine the optimal distribution of the
characteristic function
which minimises the compliance (work by external force) by iterative
update of the material distribution based on a 
sensitivity analysis.
The topology optimisation, however, suffered from numerical
instability such as mesh-dependency of optimal
configurations or checkerboard patterns.
To avoid these problems, some relaxation methods such as the density
method \cite{Bendsoe1989} or the
homogenisation method \cite{Bendsoe1988} are proposed.
Now, the topology optimisations are widely extended to design problems
in engineering fields other than structural mechanics. In the following,
we briefly review some contributions related to applications of the
topology optimisations to cloaking devices.
Andkj\ae r and Sigmund achieved to obtain some symmetric
cloak designs with dielectric elements which work for some angles of
incidence by using the density-based topology optimisation 
\cite{Andkjr2011}.
Also, Andkj\ae r et al. succeeded in designing a cloaking device that is
effective for both transverse electric (TE) and transverse magnetic (TM)
polarisations with the density-based topology optimisation \cite{Andkjr2012}.
Relaxation methods, however, allows the existence of grayscale;
region which has intermediate density between void and material
domain, which makes it difficult to fabricate the obtained
configuration. With the
help of SIMP method \cite{Bendsoe1989}, the intermediate density can be
removed but parameters for the SIMP method are chosen by a trial and
error process. With a careless choice of the parameter, the optimality
may be lost.
As one of the most promissing grayscale-free topology optimisations,
we can mention the level-set-based topology optimisation \cite{Sethian2000}.
In the level set method, boundary of design object can be expressed
clearly as a zero-level contour of a level-set
function. Also, by the use of the reaction-diffusion equation for
the update of the level-set function, complexity of the design object
can easily be adjusted \cite{Yamada2010, Choi2011}.

Many researches on the level-set based topology optimisation of cloaking
devices have been reported. 
Fujii et al. succeeded to design a cloaking device which works in
TM polarised field with the level-set based topology optimisation 
\cite{Fujii2013}.
Otomori et al. designed a cloaking device made of a ferrite
material by using the level set-based topology optimisation and
succeeded in designing a cloaking device in TM polarised field
\cite{Otomori2013}.
Their cloaking devices are, however, designed in a way to hide
a circle-shaped PEC and do not work when the shape and/or material
of the target objects are perturbed.
For engineering applications,
it is necessary to design a cloaking device which 
works stably regardless of the property of target object.
One of the most significant contribution of this paper is
to propose a new level-set-based topology optimisation for
cloaking devices which is robust to the perturbation of target
object. To this end, we modify the definition of
the objective function in a way to suppress not only scattered field
around a cloaking device but also the intensity of electro-magnetic
field in a domain where the target objects to be allocated.
By this definition, designed cloaking
device is expected to work independently on property of
target objects.

Also, previous researches on topology optimisation of cloaking devices
have a problem with respect to treatment of a design sensitivity.
In the topology optimisation, configuration of a design object is
updated repeatedly based on a topological derivative; sensitivity
of an objective function with respect to creation of an infinitesimal
circular-shaped material in a design domain.
By the use of the adjoint variable method, the topological derivative is
expressed with solutions of two boundary value problems; so called 
forward and adjoint problem. 
The performance of the topology optimisation
is highly dependent on how accurately and efficiently forward and adjoint
problem are solved.
In almost all of researches on topology optimisation, the finite
element method (FEM) is employed for the sensitivity analysis.
In the FEM, however, the infinite domain where the forward and adjoint
problem for cloaking design are often defined should be approximated
by a large but finite domain with appropriate boundary conditions.
Hence, for enough accuracy we need huge analysis domain, which
leads large numerical cost to generate a mesh and finite element
analysis.
In this study, we employ the boundary element method (BEM)
instead of the FEM for the sensitivity analysis.
In the BEM, mesh generation
is required only on the boundary. Also, the infinity itself is treated as
a boundary in the BEM and the condition at the infinity can be satisfied
strictly with the help of the Green function.
Our previous researches have shown the effectiveness of the level-set
based topology optimisation with the BEM in the electromagnetic field
\cite{Isakari2016multi}, 
heat conduction problem \cite{Jing2013, Jing2015} and sound field problem
\cite{Isakari2014, Isakari2016level}.
We further improve our previous BEM-based optimisations to
develop an efficient computation of the
topological derivative.
The forward and adjoint problem can be reduced to algebraic
equations by a numerical analysis method such as the FEM or BEM.
These two equations can be solved efficiently with direct solvers such
as the LU decomposition
since the coefficient matrices of these equations are the same
in many cases.
The coefficient matrix derived by the BEM is, however, fully-populated,
which causes large numerical cost for generation of the coefficient
matrix and the LU decomposition. 
Hence for the fast computation of the topological derivative, 
acceleration of the BEM and LU decomposition is necessary.
The fast multipole method (FMM) is well known acceleration method for
the BEM \cite{Greengard1987fast, Rokhlin1990rapid}.
The FMM is usually combined with iterative solvers for algebraic 
equations because it is a method to accelerate the matrix-vector
product of the BEM coefficient matrix with an arbitrary vector.
Hence, we need to solve the forward and adjoint problem 
individually. Also, the number of iterations to obtain the solution
is highly dependent on the property of problems such as complexity
of geometry and material constant. 
In this study, as an alternative acceleration method for the BEM, we
employ the ${\mathcal H}$-matrix method 
\cite{Hackbusch1999sparse, Hackbusch2000479}
which is based on
hierarchical blocking of the coefficient matrix and low rank
approximation.
This method can easily be
combined with direct solvers.
Furthermore, we can reduce the cost to generate the coefficient matrix
and memory to store the matrix by using the adaptive cross approximation
(ACA) for the low rank approximation.
Hence, it is expected that we can accelerate the sensitivity analysis
by using the efficient computation of the coefficient matrix with the
ACA and the accelerated LU decomposition by the ${\mathcal H}$-matrix method.

The rest of this paper is organised as follows: In the second section,
we formulate the level-set based topology optimisation method for
cloaking devices. In the third section, we show the derivation of the
topological derivative. In the fourth section, we explain about the 
electro-magnetic field analysis with the BEM.  In the fifth
section, we introduce an efficient computation of the topological
derivative with the ${\mathcal H}$-matrix method. In the sixth section,
we show the effectiveness of proposed methods with some numerical
examples.

\section{Topology optimisation of cloaking devices}
In this section, we define conventional and proposed optimisation
problems to design cloaking devices.
We also show the procedure for solving the optimisation problems
including configuration expression and its update.

\subsection{A conventional optimisation problem in design for cloaking
  devices and its modification}
\begin{figure}[htbp]
  \begin{center}
    \includegraphics[scale=0.2]{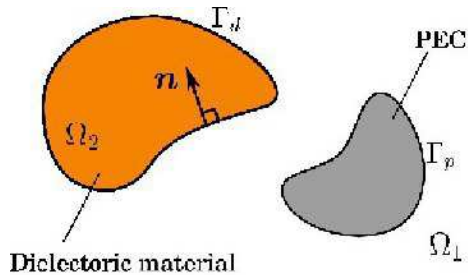}
    \caption{Definition of symbols.}
    \label{fig.die_pec}
  \end{center}
\end{figure}
\begin{figure}[!h]
 \begin{center}
  \includegraphics[scale=0.18]{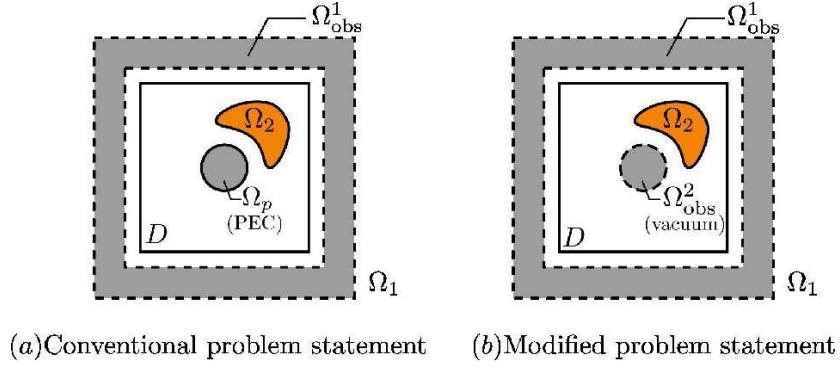}
  \caption{Two types of optimisation problems in design for cloaking devices.}
  \label{fig.opt_cloak}
 \end{center}
\end{figure}

In this study, we consider optimisation problems of cloaking 
devices in two-dimensional TM polarised
field and assume that cloaking device is made of dielectric
elements.
For the time being, we assume that objects to be hidden is 
made of perfectly electric conductor (PEC). 
This restriction will be removed in the proposed optimisation method.
We denote a vacuum domain, dielectric elements
and PEC as
$\Omega_1$, $\Omega_2$ and $\Omega_p$, respectively. 
Also, we express boundaries of $\Omega_2$ and $\Omega_p$ by $\Gamma_d$
, $\Gamma_p$ (\figref{fig.die_pec}).
Purpose of cloaking devices is to make the PEC invisible by 
suppressing the scattered field around the PEC.
Previous researches have achieved this by determining
configuration of cloaking devices which reduce the intensity of 
scattered field.
In this case, the design problem can be formulated as the following
optimisation problem:
\begin{alignat}{2}
 &\mathrm{Find}\ \Omega_2\ \mathrm{in}\ D \notag \\
 &\mathrm{such\ that}\min J \notag \\
 &J=\sum_{m=1}^M\| u({\bs x}_m^1) - u^\mathrm{inc}({\bs
 x}_m^1)\|^2 \label{eq.conv_obj}\\
 &\mathrm{subject\ to} \\
 &\nabla^2 u(\bs x) + k_1^2 u(\bs x) = 0 & &\ \ \ \bs x \in \Omega_1 \label{tmbvp1} ,\\
 &\nabla^2 u(\bs x) + k_2^2 u(\bs x) = 0  & &\ \ \  \bs x \in \Omega_2 \label{tmbvp2} ,\\
 &u^1=u^2 & &\ \ \  \bs x \in \Gamma_d \label{tmbvp3} ,\\
 &\frac{1}{\mu_1}\left(\frac{\partial u}{\partial n}\right)^1
 =\frac{1}{\mu_2}\left(\frac{\partial u}{\partial n}\right)^2
 & &\ \ \  \bs x \in \Gamma_d \label{tmbvp4} ,\\
 &u = 0 & &\ \ \  \bs x \in \Gamma_p \label{tmbvp5} ,\\
 &\mathrm{Radiation\ condition} & &\ \ \  \bs x \rightarrow \infty ,
 \label{tmbvp6}
\end{alignat}
where $u$ (resp. $u^\mathrm{inc}$) denotes the total electric field
(resp. incident field), and $D$ denotes a fixed design domain in which a
circle-shaped PEC $\Omega_p$ is allocated (\figref{fig.opt_cloak}, left).
Also, ${\bs x}_m^1\ (m=1,\cdots,M)$ denotes an observation point
allocated in an observation domain $\Omega_\mathrm{obs}^1$ around the
design domain $D$, and $M$ is the number of the observation points.
In the boundary conditions (\ref{tmbvp3}) and
(\ref{tmbvp4}) on $\Gamma_d$,
the variable with the superscripted index $i$ denotes the limit of the
variable to the boundary $\Gamma_d$ from $\Omega_i$.
Also, $\mu_i$ denotes the permeability of domain $\Omega_i$, and
$\partial/\partial n := \bs n\cdot \nabla$ 
denotes an outward normal derivative from $\Omega_1$.
$k_i$ in equations (\ref{tmbvp1}) and (\ref{tmbvp2})
is the wave number defined as follows with angular
frequency $\omega$ and permittivity $\varepsilon_i$ in $\Omega_i$:
\begin{equation}
 k_i=\omega \sqrt{\mu_i\varepsilon_i}.
\end{equation}

The optimal configuration of this problem is, however, highly
specialised to the circular PEC, and the effectiveness
is declined when the PEC is replaced by other shaped one.
To avoid this problem, we modify the optimisation
problem. 
We augment the objective function by adding
the sum of the intensity of the electric field in a fixed
circular domain $\Omega_\mathrm{obs}^2$ which covers
the PEC $\Omega_p$ in the conventional optimisation problem.
i.e., the optimisation problem is defined as follows:
\begin{alignat}{2}
 &\mathrm{Find}\ \Omega_2\ \mathrm{in}\ D \notag \\
 &\mathrm{such\ that}\min J \notag \\
 &J=
 \sum_{m=1}^M\| u({\bs x}_m^1) - u^\mathrm{inc}
 ({\bs x}_m^1)\|^2
 +
 \sum_{n=1}^N\| u({\bs x}_n^2) \|^2 \label{mod_obj} \\
 &\mathrm{subject\ to},\\
 &\nabla^2 u(\bs x) + k_1^2 u(\bs x) = 0 & &\ \ \ \bs x \in \Omega_1
 \cup \Omega_\mathrm{obs}^2, \label{mod_bvp1}\\
 &\nabla^2 u(\bs x) + k_2^2 u(\bs x) = 0  & &\ \ \  \bs x \in \Omega_2 ,\\
 &u^1=u^2 & &\ \ \  \bs x \in \Gamma_d,\\
 &\frac{1}{\mu_1}\left(\frac{\partial u}{\partial n}\right)^1
 =\frac{1}{\mu_2}\left(\frac{\partial u}{\partial n}\right)^2 
 & &\ \ \  \bs x \in \Gamma_d ,\\
 &\mathrm{Radiation\ condition} & &\ \ \  \bs x \rightarrow \infty \label{mod_bvp2},
\end{alignat}
where $\bs x_n^2\ (n=1,\cdots,N)$ denotes an observation point in 
$\Omega_\mathrm{obs}^2$, and $N$ is the number of observation points.
By the above definition, it is expected to obtain the cloaking device
which works regardless of the PEC shape since scattering on the
surface of PEC allocated in $\Omega_\mathrm{obs}^1$ will be almost
zero. Furthermore, it is expected that the obtained cloaking device
work even when the PEC is replaced by dielectric one.

In the following subsections, a procedure to solve the optimisation
problems is presented.

\subsection{Expression of configuration with the level-set method}
\begin{figure}[!h]
 \begin{center}
  \includegraphics[scale=0.13]{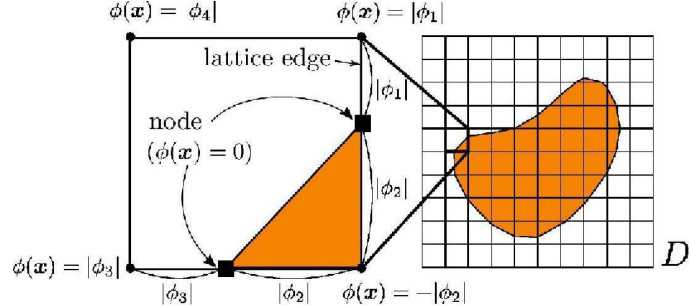}
  \caption{Generation of the boundary mesh.}
  \label{fig.lsf2bndry}
 \end{center}
\end{figure}

We use the level-set method to express the configuration
of a design object $\Omega_2$. In the level-set method, boundary of a design object
$\Gamma_d$ is expressed as zero-level contour of a level-set function
$\phi(\bs x)$ defined in the design domain. Inner and outer domain
of the design object
is distinguished by the value of the level set function as follows:
\begin{alignat}{2}
 &\Omega_1 & &= \{\ \bs x \ |\   0 < \phi(\bs x) \le 1\ \}  , \\
 &\Gamma_d & &= \{\ \bs x \ |\  \phi(\bs x)=0\  \}          ,  \\
 &\Omega_2 & &= \{\ \bs x \ |\  -1 \le \phi(\bs x) < 0\  \} .
\end{alignat}
The value of $\phi(\bs x)$ is stored on a grid point of lattice
expanding the design domain $D$.
In order to generate a boundary mesh, we initially interpolate the value
of the level-set function between neighbouring lattice points by a line.
By connecting points at which $\phi(\bs x)=0$ in each lattice, boundary
elements are generated (\figref{fig.lsf2bndry}).
Generated boundary mesh by the above process may, however, includes
boundary elements whose length are uneven.
Hence, after the mesh generation, we improve these boundary elements so that
the length of each boundary element becomes almost the same \cite{Isakari2016multi}.

\subsection{Update of configuration with the reaction-diffusion
  equation}
  For the update of the level-set function at each optimisation step, we
  solve the following reaction-diffusion equation in the design
  domain $D$ with appropriate initial and boundary conditions with
  respect to ficticious time $t$ \cite{Yamada2010}.
 \begin{align}
  \pd{\phi}{t} =C \textcolor{black}{\mathrm{sgn}(\phi(\bs x))}{\mathcal
  T}(\bs x) +   \tau l^2 \nabla^2 \phi(\bs x). \label{rdeq}
 \end{align}
 In the right hand side (RHS) of equation (\ref{rdeq}), the first term
 denotes the
 direction and scale of the update of $\phi(\bs x)$. $C$ is a constant
 which arranges the scale. ${\mathcal T}(\bs x)$ denotes a topological
 derivative which is a sensitivity
 of the objective function with respect to new creation of an
 infinitesimal circular dielectric or vacuum domain. 
 For the derivation of the topological derivative,
 the reader is referred to Section $3$.
 The second term of the RHS in equation (\ref{rdeq}) is so called
 Tikhonov's regularisation term and 
 works as a perimeter constraint for a design object. 
 By adjusting a parameter $\tau$,
 the complexity of obtained configuration after the update can
 be arranged \cite{Choi2011}. $l$ is a characteristic length, in this
 study the side length of $D$ is employed.

\subsection{Algorithm of the topology optimisation} 
The algorithm of our topology optimisation is summerised as follows:
\begin{itemize}
 \item Step 1: Initialise a level-set function in accordance with a
       initial configuration.
 \item Step 2: Generate a boundary mesh based on the level-set function
       (see Section $2.2$).
 \item Step 3: Solve the boundary value problem 
       (either (\ref{tmbvp1})--(\ref{tmbvp6}) or
       (\ref{mod_bvp1})--(\ref{mod_bvp2})) with the BEM
       (see Section $4$).
 \item Step 4: Compute the objective function and check if the convergence
       condition is satisfied or not. In this study we judge the
       objective function is
       converged at step $k$ when the objective function of the last $50$
       steps $(J_{k-49},\cdots,J_{k})$ satisfy the following relation:
       \begin{align}
	&\ \ J' \le \varepsilon_1^\mathrm{conv},\\
	&\max_{1\le i,j \le 50}
	\left|
	\frac{J_{k-i+1}}{J_{k-j+1}}
	\right|
	\le\varepsilon_2^\mathrm{conv},
       \end{align}
       where $J'$ denotes the gradient of the objective function which
       is evaluated with the least square for 
       $\log (J_{k-49}),\cdots, \log (J_k)$. Also,
       $\varepsilon_1^\mathrm{conv}$, $\varepsilon_2^\mathrm{conv}$ are
       parameters.
       When the condition is satisfied, we terminate the 
       optimisation process.
 \item Step 5: Compute the topological derivative derived in the Section $4$.
 \item Step 6: Update the distribution of the level-set function by solving
       the reaction and diffusion equation (\ref{rdeq}) with the FEM
       (see Section 2.3) and back to the step 2.
\end{itemize}
\section{Topological derivative}
\begin{figure}[htbp]
  \begin{center}
    \includegraphics[scale=0.12]{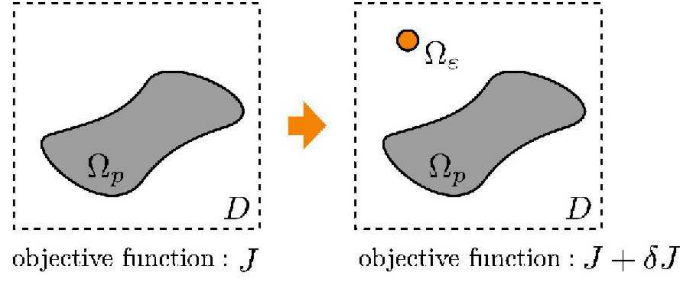}
    \caption{Allocation of an infinitesimal circle in the design domain.}
    \label{fig.mini_circle}
  \end{center}
\end{figure}

In our topology optimisation, we update a configuration of
a design object
iteratively based on a topological derivative
(see the reaction-diffusion equation (\ref{rdeq})); sensitivity of an
objective function $J$ when an infinitesimal circle
$\Omega_\varepsilon$ is allocated in a
design domain $D$ (\figref{fig.mini_circle}) \cite{Carpio2008}.
In this section, we consider an objective function defined as summation of a
function of $u$ as follows:
\begin{align}
 J=\sum_{m=1}^M f(u(\bs x_m^\mathrm{obs})) ,
\end{align}
where $\bs x_m^\mathrm{obs} \not\in D\ (m=1,\cdots, M)$ is an
observation points.
Other types of the objective function such as the one defined on the
boundary can readily be obtained with a similar procedure as the
following discussion. 

Topological derivative is defined as a first coefficient in 
asymptotic expansion of the objective function by measure of the
allocated infinitesimal circle $s(\varepsilon)$ whose centre is $\bs x$:
\begin{equation}
 \delta J 
  =
  {\mathcal T}({\bs x}) s(\varepsilon)
  +
  o(s(\varepsilon)) \label{tddef} ,
\end{equation}
in which $\delta J$ is a variation of the objective function which 
is evaluated as follows:
\begin{align}
 \delta J
  &=
  \frac{\partial J}{\partial u_r}\delta u_r
  +
  \frac{\partial J}{\partial u_i}\delta u_i \notag \\
  &=
  \Re \left[\frac{\partial J}{\partial u}\delta u\right] \label{dj} ,
\end{align}
where $\delta u$ is a variation of the electric response when
the infinitesimal circular domain is allocated on $\bs x$. Also,
indices `$r$' and `$i$'
represent the real part and imaginary part of $\delta u$, respectively.
In the case that the allocated infinitesimal circle is made
of dielectric element, 
it can easily be confirmed that $\delta u$ satisfies the following boundary
value problem:
\begin{alignat}{2}
 &\nabla^2 \delta u({\bs x})+k_1^2 \delta u({\bs x}) 
 = 0
 &{\bs x} &\in \Omega_1 ,    \label{bvp_du1}\\
 &\nabla^2 \delta u({\bs x})+k_2^2 \delta u({\bs x}) 
 = 0
 &{\bs x} &\in \Omega_2 ,   \\
 &\nabla^2 \hat u(\bs x) + k_2^2  \hat u (\bs x)
 = 0
 &{\bs x} &\in \Omega_\varepsilon, \\ 
 &\delta u_1 = \delta u_2
 &{\bs x} &\in \Gamma_d, \\
 &\frac{1}{\mu_1}
 \left(
 \frac{\partial \delta u}{\partial n}
 \right)^1
 =
 \frac{1}{\mu_2}
 \left(
 \frac{\partial \delta u}{\partial n}
 \right)^2
 &{\bs x} &\in \Gamma_d, \\
 &\delta u = 0
 &{\bs x} &\in \Gamma_p, \\
 &u_1 + \delta u_1 = \hat u
 &{\bs x} &\in \Gamma_\varepsilon, \label{uhat1}\\
 &\frac{1}{\mu_1}
 \frac{\partial (u_1+\delta u_1)}{\partial n}
 =
 \frac{1}{\mu_2}
 \frac{\partial \hat u}{\partial n}
 &{\bs x} &\in \Gamma_\varepsilon \label{uhat2} ,\\
 &\mathrm{Radiation\ condition}
 &|\bs x| &\rightarrow \infty , \label{bvp_du3}
\end{alignat}
where $\hat u$ is the electro-magnetic response in $\Omega_\varepsilon$.
We define the normal vector $\bs n$ is positive when $\bs n$ is directed
from $\Omega_1$.
In this study, we employ the adjoint variable method to evaluate the
RHS of equation (\ref{dj}). Namely, we define the following
adjoint problem:
\begin{alignat}{2}
 &\nabla^2 \tilde u(\bs x) + k_1^2 \tilde u(\bs x) +
 \sum_{m=1}^M
 \frac{\partial f(u(\bs x_m^\mathrm{obs}))}
 {\partial u}
 \delta(\bs x-\bs x^\mathrm{obs}_m)
 = 0
 &\ \ \ \bs x &\in \Omega_1,  \label{adj1} \\
 &\nabla^2 \tilde u(\bs x) + k_2^2 \tilde u(\bs x) = 0
 &\ \ \ \bs x &\in \Omega_2, \\
 &\tilde u^1 = \tilde u^2
 &\ \ \ \bs x &\in \Gamma_d,  \\
 &\frac{1}{\mu_1}
 \left(
 \frac{\partial \tilde u}{\partial n}
 \right)^1
 =
 \frac{\textcolor{black}{1}}{\mu_2}
 \left(
 \frac{\partial \tilde u}{\partial n}
 \right)^2
 &\ \ \ \bs x &\in \Gamma_d, \\
 &\tilde u = 0
 &\ \ \ \bs x &\in \Gamma_p,  \\
 &\mathrm{Radiation\ condition} 
 &\ \ \ |\bs x| &\rightarrow \infty. \label{adj2}
\end{alignat}
With the help of the reciprocal theorem for
\begin{itemize}
 \item $\tilde u$ and $\delta u$ in 
       $\Omega_1\backslash \overline{\Omega_\varepsilon}$,
 \item $\tilde u$ and $\delta u$ in $\Omega_2$,
 \item $\tilde u$ and $\hat u$ in $\Omega_\varepsilon$,
\end{itemize}
the RHS of equation (\ref{dj}) is evaluated as follows:
\begin{align}
 \delta J
 =
 \Re\left[
 \int_{\Omega_\varepsilon}
 \omega^2(\varepsilon_2-\varepsilon_1)
 \hat{u}\tilde{u}
 ~\mathrm{d}\Omega 
 \right] \label{deltaj_exp}.
\end{align}
The asymptotic behaviour of $\hat u$ and $\tilde u$ in
$\Omega_\varepsilon$ is evaluated as follows \cite{Carpio2008}:
\begin{alignat}{2}
 &\hat u(\bs x)
 =
 u(\bs x_0)
 +
 \frac{2\varepsilon_2}{\varepsilon_1+\varepsilon_2}
 u_{,j}(\bs x_0)(\bs x - \bs x_0)_j
 +
 o(\varepsilon)
 & &\ \ \ 
 \bs x \in \Omega_\varepsilon, \label{uhat_exp} \\
 &\tilde u(\bs x)
 =
 \tilde u(\bs x_0)
 +
 \tilde u_{,j}(\bs x_0)(\bs x - \bs x_0)_j
 +
 o(\varepsilon)
 & &\ \ \ 
 \bs x \in \Omega_\varepsilon. \label{tilde_exp}
\end{alignat}
By substituting (\ref{uhat_exp}), (\ref{tilde_exp}) into
(\ref{deltaj_exp}), we obtain the following expression.
\begin{align}
 \delta J
 =
 \Re\left[
 \left(
 \omega^2 (\varepsilon_2-\varepsilon_1)
 \hat{u}\tilde{u}
 \right)
 \pi\varepsilon^2
 +
 o(\varepsilon^2)
 \right].
\end{align}
From the definition (\ref{tddef}), the topological derivative is
derived as follows:
\begin{equation}
 {\mathcal T}(\bs x)=\Re\
  \left[
   \omega^2 (\varepsilon_2-\varepsilon_1)
   u(\bs x) \tilde{u}(\bs x)
  \right]. \label{tdv}
\end{equation}

\section{Electromagnetic field analysis in the 2D infinite domain}
In this section, we show a method to solve two dimensional Maxwell's equations in the
infinite domain which appears in the constraint in the optimisation
problem (\ref{tmbvp1})--(\ref{tmbvp6}) with the boundary element method (BEM). 
The BVP (\ref{mod_bvp1})--(\ref{mod_bvp2}) can similarly be solved.
The solution of the BVP (\ref{tmbvp1})--(\ref{tmbvp6}) has the
following integral representation:
\begin{alignat}{2}
u(\bs x)
&=
u^\mathrm {inc}(\bs x)
-
\int_{\Gamma_d}
\frac{\partial G^1(\bs x,  \bs y)}{\partial n_y}
u(\bs y) ~\mathrm{d}\Gamma_y 
+\int_{\Gamma_p \cup \Gamma_d} G^1(\bs x,  \bs y)
 \frac{\partial u(\bs y)}{\partial n_y} ~\mathrm{d}\Gamma_y
 &&\ \ \ \bs x \in \Omega_1 \label{inner1} ,\\
u(\bs x)
&=
\int_{\Gamma_d}
\frac{\partial G^2(\bs x,  \bs y)}{\partial n_y}
u(\bs y) ~\mathrm{d}\Gamma_y 
-\int_{\Gamma_d} G^2(\bs x,  \bs y)
 \frac{\partial u(\bs y)}{\partial n_y} ~\mathrm{d}\Gamma_y
 &&\ \ \ \bs x \in \Omega_2 \label{inner2},
\end{alignat}
where $G^i(\bs x, \bs y)\ (i=1, 2)$ denotes the fundamental solution of the Helmholtz equation
in 2D which is expressed by the Hankel function of the first kind of order 0 as follows:
\begin{align}
 G^i(\bs x, \bs y)=\frac{i}{4}H_0^{(1)}(k_i|\bs x - \bs y|).
\end{align}
By taking the limit as $\bs x \rightarrow \Gamma_p$ and $\Gamma_d$, we
can solve the obtained boundary integral equations.
The obtained equations may, however, suffer from so called fictitious
eigenfrequency problem.
In order to avoid the problem, we employ the PMCHWT formulation
\cite{Chew1995} and obtain the following boundary integral equations:
\renewcommand{\arraystretch}{2.3}
\begin{equation}
 \left[
  \begin{array}{ccc}
  \displaystyle
   -{\mathcal S}^{\Gamma_p}_1           & 
   \displaystyle \frac{1}{\mu_1}{\mathcal D}^{\Gamma_d}_1   &
   -\mu_1 {\mathcal S}^{\Gamma_d}_1 \\

  \displaystyle
   -{\mathcal D}_1^{\Gamma_p} &
   \displaystyle
   \frac{{\mathcal N}_1^{\Gamma_d}}{\mu_1}+\frac{{\mathcal N}_2^{\Gamma_d}}{\mu_2} &
   -({\mathcal D}_{1}^{\Gamma_d*}+{\mathcal D}_{2}^{\Gamma_d*})\\
   -\mu_1{\mathcal S}^{\Gamma_p}_1 & {\mathcal D}^{\Gamma_d}_1+{\mathcal D}^{\Gamma_d}_2 &
   -(\mu_1{\mathcal S}_1^{\Gamma_d}
   +\mu_2 {\mathcal S}_2^{\Gamma_d} )
  \end{array}
 \right]
 \left[
 \begin{array}{c}
  w_p \\
  u_d \\
  w_d 
 \end{array}
 \right]
 =
 \left[
  \begin{array}{c}
   w_p^\mathrm{inc} \\
   w_d^\mathrm{inc} \\
   u_d^\mathrm{inc} 
  \end{array}
 \right] \label{eq.pmchwt_tm} ,
\end{equation}
where indices `$p$' and `$d$' indicate the variable defined on $\Gamma_p$
and $\Gamma_d$, respectively. 
Also, $w_p$, $w_p^\mathrm{inc}$, $w_d$ and $w_d^\mathrm{inc}$ are
defined as follows:
\begin{align}
 &w_p
 =
 \frac{1}{\mu_1}
 \left(
 \frac{\partial u}{\partial n}
 \right)^1, \\
 &w_p^\mathrm{inc}
 =
 \frac{1}{\mu_1}
 u_p^\mathrm{inc}, \\
 &w_d
 =\frac{1}{\mu_1}
 \left(
 \frac{\partial u}{\partial n}
 \right)^1
 =\frac{1}{\mu_2}
 \left(
 \frac{\partial u}{\partial n}
 \right)^2 , \\
 &w_d^\mathrm{inc}
 =\frac{1}{\mu_1}
 \frac{\partial u^\mathrm{inc}}{\partial n} .
\end{align}
${\mathcal D}_i^{\Gamma},\ {\mathcal S}_i^{\Gamma},\ {\mathcal D}_{i}^{\Gamma *}$ 
and $\ {\mathcal N}_i^{\Gamma}$ respectively denote the following operators:
\begin{alignat}{2}
 &[{\mathcal D}_i^{\Gamma} \phi](\bs x)& &=\int_\Gamma
 \frac{\partial G^i(\bs x, \bs y)}{\partial n_y}
 \phi(\bs y) ~\mathrm{d}\Gamma_y ,\label{eq:sss}\\
 &[{\mathcal S}_i^{\Gamma} \psi](\bs x)& &=\int_\Gamma
 G^i(\bs x,\bs y)
 \psi(\bs y)
 ~\mathrm{d}\Gamma_y ,\\
 &[{\mathcal D}_{i}^{\Gamma *} \psi](\bs x)& &=\int_\Gamma
 \frac{\partial G^i(\bs x,\bs y)}{\partial n_x}
 \psi(\bs y)
 ~\mathrm{d}\Gamma_y ,\\
 &[{\mathcal N}_i^{\Gamma} \phi](\bs x)& &=\int_\Gamma
 \frac{\partial^2 G^i(\bs x,\bs y)}{\partial n_x
 \partial n_y}
 \phi(\bs y) ~\mathrm{d}\Gamma_y,\label{eq:nnn}
\end{alignat}
where $\phi$ and $\psi$ are density functions.
By discretising the boundary $\Gamma_p \cup \Gamma_d$ and the electric
field and magnetic field $(u, w)$ in the boundary integral equation
(\ref{eq.pmchwt_tm}) by
linear and constant element, respectively, the collocation gives
the algebraic equations. After computing the electro-magnetic field on
the boundary by solving the algebraic equations, we obtain the
electric response in $\Omega_1$ and $\Omega_2$ by substituting the
solutions on the boundary to equations (\ref{inner1}) and (\ref{inner2})
, respectively.

\section{Efficient sensitivity analysis with the $\mathcal H$-matrix method}
In the previous section, we showed the expression of the topological
derivative with the adjoint variable method.
For the computation of the topological derivative, we need to
compute the electro-magnetic responses for two incident field; forward
and adjoint.
For the efficient computation of the topological derivative, 
fast direct solver is suitable since the coefficient matrix of
forward and adjoint problem is the same, and we can solve
these problems at the same time with a direct solver.
The hierarchical matrix method; so called ${\mathcal H}$-matrix method 
is one of the promissing acceleration methods for matrix operations
\cite{Hackbusch1999sparse, Hackbusch2000479}, and we can solve
the algebraic equations with $O(N\log N)$ computational cost 
by using accelerated LU decomposition with the ${\mathcal H}$-matrix
method (${\mathcal H}$LU) in which $N$ denotes the dgree of freedom. 

In the ${\mathcal H}$-matrix method, we firstly express a matrix in the
form of the ${\mathcal H}$-matrix which is a matrix constructed through
the following two steps:
\begin{itemize}
 \item Hierarchical blocking of a matrix.
 \item Low rank approximation of submatrices which express influence from
       far field.
\end{itemize}
By applying efficient operations to each sub low rank matrices,
the computational cost for matrix operations is reduced.
In this section, we explain the construction of the ${\mathcal H}$
-matrix following the above two steps.
Also, we show an efficient sensitivity analysis method 
with the ${\mathcal H}$-matrix method.

\subsection{Blocking of a coefficient matrix}
\begin{figure}[htbp]
  \begin{center}
   \includegraphics[scale=0.4]{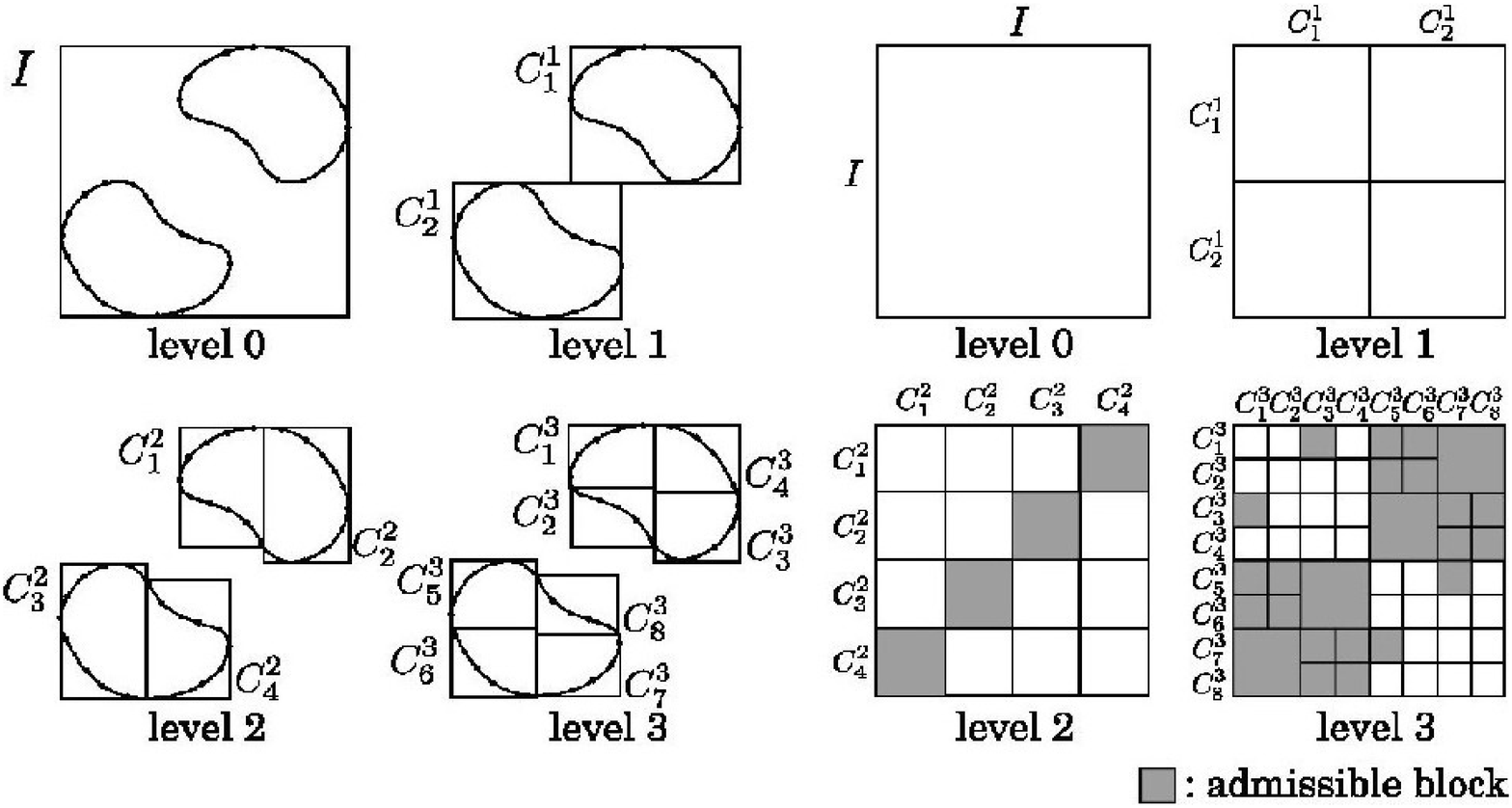}
   \caption{An example of boundary clustering and corresponding
   division of a coefficient matrix.}
   \label{fig.cluster}
  \end{center}
\end{figure}

\begin{figure}[htbp]
  \begin{center}
    \includegraphics[scale=0.18]{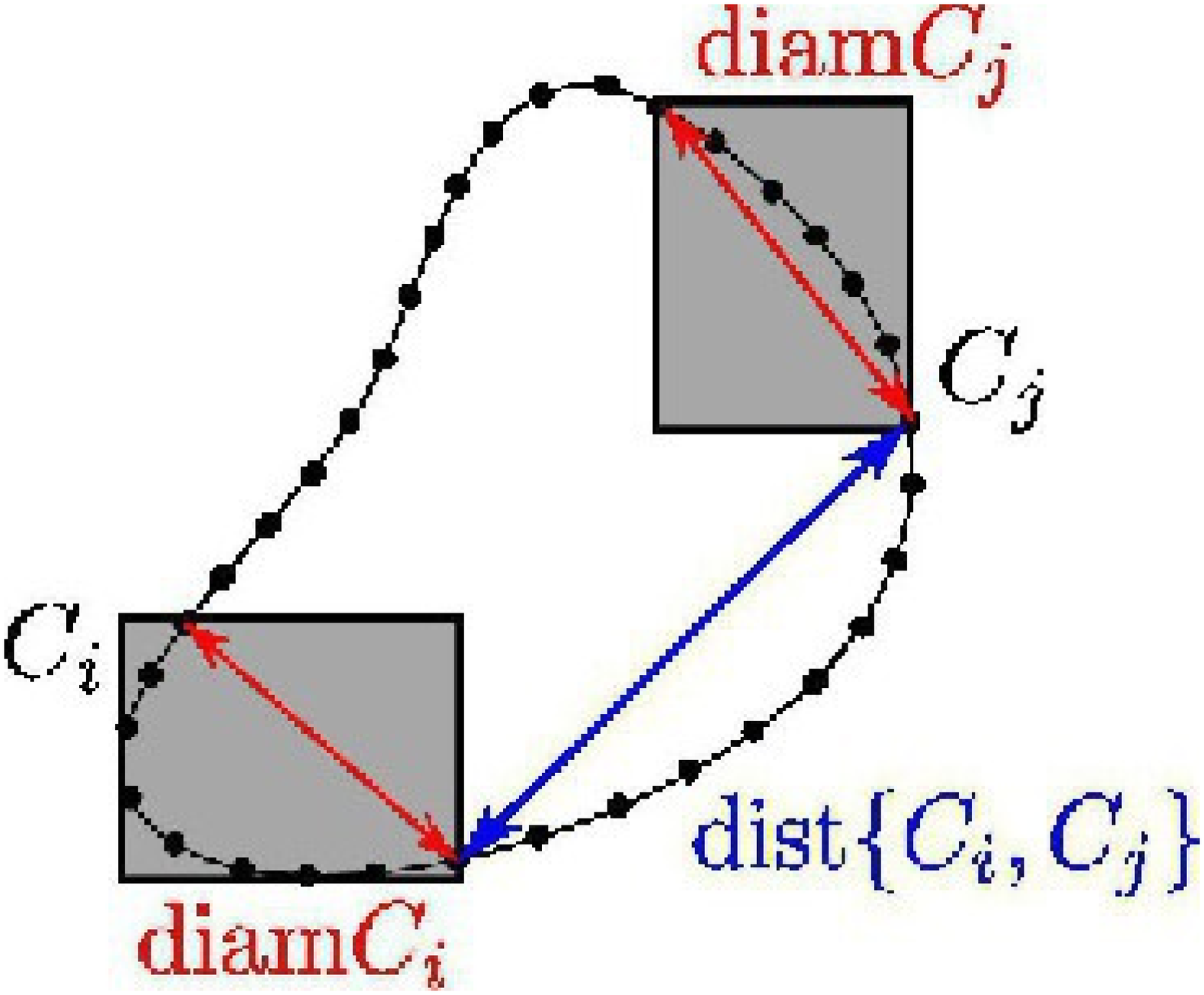}
    \caption{Definition of $\mathrm{diam}$ and $\mathrm{dist}$.}
    \label{fig.admis}
  \end{center}
\end{figure}

As an example, we consider a vacuum domain $\Omega_1$ in which some PEC
objects ${\Omega_p}$ are allocated. The boundary integral equation for
this problem is derived as follows:
\begin{align}
u^\mathrm {inc}(\bs x)
&=
-\int_{\Gamma_p}
G^1(\bs x,  \bs y)
\frac{\partial u(\bs y)}{\partial n_y}
 ~\mathrm{d}\Gamma_y
 \ \ \ \bs x \in \Gamma_p \label{bie_eg}. 
\end{align}
By discretising the boundary and $\partial u/ \partial n_y$ of
(\ref{bie_eg}) respectively with linear and constant element, we obtain
a dense coefficient matrix $A$ whose entry $a_{ij}$ is expressed as
follows:
\begin{align}
 a_{ij}=-\int_{\Gamma_j} 
 G^1(\bs x_i, \bs y)
 ~\mathrm{d}\Gamma_y, \label{aij}
\end{align}
where $\Gamma_j$ denotes the $j$th discretised boundary element,
and $\bs x_i$ is the collocation point on $\Gamma_i$.
One finds that each entry $a_{ij}$ is associated with a collocation $\bs x_i$
and a boundary element $\Gamma_j$.
Hence, blocking of $A$ can be
done based on partition of $\Gamma_p$.
Namely, we generate a boundary cluster and divide the coefficient matrix
based on the clustering by the following strategy:
\begin{itemize}
 \item Make a rectangle enclosing a boundary which we are going to be divided.
       We define a set of elements in the rectangle as a `cluster'. Also, 
       we denote the number of division which is required to obtain the 
       cluster as a `level' to which the cluster belongs.
       In the following, we denote $i$th cluster in the level $l$ as $C_i^l$.
 \item Divide the longer side of the rectangle into two parts. We define 
       sets of elements in each obtained rectangle as new clusters
       (\figref{fig.cluster}). The coefficient matrix
       is also divided corresponding to the division of the boundary.
       Each block matrix corresponds to a combination of two clusters
       in the same level.
 \item When two clusters $C_i^l$ and $C_j^l$ at the same level $l$
       satisfy the following
       admissibility condition, we define the block matrix 
       $C_i^l\times C_j^l$
       as an admissible block (the corresponding submatrix represents
       far-field influence) and stop dividing the
       block (\figref{fig.cluster}).
       \begin{align}
	\mathrm{min}\{\mathrm{diam}\ C_i^l, \mathrm{diam}\ C_j^l\}
	\le
	\eta\ \mathrm{dist}\{C_i^l, C_j^l\} \label{eq.admis},
       \end{align}
       where $\eta$ denotes a real constant which sets the strictness of
       the condition (The larger $\eta$ is set, the more blocks are
       recognised  as admissible). `$\mathrm{diam}$' and `$\mathrm{dist}$'
       respectively denote diameter of a cluster and distance of two
       clusters (\figref{fig.admis})
       and defined as follows:
       \begin{align}
	&\mathrm{diam}\ C_i^l = \max_{\bs x, \bs y \in C_i^l}|\bs x-\bs y|, \\
	&\mathrm{dist}\{C_i^l,C_j^l\} = \min_{\bs x \in C_i^l, \bs y \in
	C_j^l} |\bs x-\bs y|.
       \end{align}
 \item We repeat the above process until the number of nodes in each
       cluster is less than a preset parameter $n_\mathrm{min}$.
       We define such clusters as `leaf cluster'. When at least
       one of two clusters $C_s^l$, $C_t^l$ at the same level $l$ is a
       leaf cluster and $C_s^l$, $C_t^l$ does
       not satisfy the admissibility condition (\ref{eq.admis}), 
       we define the block $C_s^l\times C_t^l$ as an inadmissible block
       (the corresponding submatrix represents near-field influence).
\end{itemize}
After we obtain the structure of the ${\mathcal H}$-matrix,
we compute the submatrices corresponding to inadmissible and admissible
blocks. While the inadmissible matrices are calculated according to its
definition (\ref{aij}), the admissible counterparts are evaluated in a
low rank approximated form which is presented in the following subsection.


\subsection{Efficient computation of block matrices with the ACA and parallelisation}
Low rank approximation is a method to express a matrix $A \in
{\mathbb C}^{m\times n}$ with two vectors $\bs a_l \in {\mathbb C}^{m\times 1}$ 
and $\bs b_l \in {\mathbb C}^{n\times 1}$ $(l=1,\cdots,k)$ as follows:
\begin{align}
 A \simeq S_k = \sum_{i=l}^k \bs a_l \bs b_l^H, \label{lowrank}
\end{align}
where $S_k$ is an approximation matrix of $A$, and $k$ is the rank of $S_k$.
Although the original matrix $A$ requires $O(mn)$ memory to store the matrix,
approximated matrix $S_k$ requires $O(k(m+n))$ memory. Hence, when the rank
$k$ satisfies $k < mn/(m+n)$ we can reduce the memory consumption.
Also, we can define efficient matrix operations such as
matrix-vector product, LU decomposition, etc for matrices expressed
in the form of (\ref{lowrank}).

The adaptive cross approximation (ACA) is one of the techniques to realise
the low-rank approximation (\ref{lowrank}) of matrices.
The algorithm of the ACA is shown in Algorithm \ref{alg_aca}. 
For more details such as the selection of the initial row and
efficient evaluation of the Frobenius norm etc, the reader may
consult \cite{bebendorf2008hierarchical}.
              
As a result of approximation with the ACA, the approximation matrix
satisfies the following error estimation:
\begin{align}
 \|A - S_k\|_F \le \varepsilon \|A\|_F ,
\end{align}
where $\|\cdot \|_F$ denotes the Frobenius norm, and $\varepsilon$ is a
tolerance for the ACA.
In the algorithm of the ACA, we need not to compute all of the 
entries of $A$, which can reduce the assemble cost of the coefficient
matrix.

\begin{algorithm}
 \caption{Adaptive cross approximation (ACA)}
 \label{alg_aca}
 \begin{algorithmic}
  \STATE Let $i_1=${\it Initial Row}
  \FOR{$k=1\ \mathrm{to}\ m$}
  \STATE $\bs b_k = A(i_k,:) - \sum_{l=1}^{k-1}(\bs a_l)_{i_k} \bs b_l$
  \STATE $j_k = \mathrm{argmax}(|(\bs b_k)_j|)$
  \STATE $\bs a_k = A(:,j_k) - \sum_{l=1}^{k-1}(\bs b_l)_{j_k} \bs a_l$
  \STATE $\bs a_k = (\bs b_k)^{-1} \bs a_k$
  \IF {$(\|\bs a_k\|_F \|\bs b_k\|_F \le \varepsilon \|A_k\|_F)$}
  \STATE exit
  \ENDIF
  \STATE $i_{k+1}=\mathrm{argmax}(|(\bs a_k)_i|)$
  \ENDFOR
 \end{algorithmic}
\end{algorithm}

For further acceleration of the ${\mathcal H}$-matrix generation process, we
consider to parallelise the procedure.
Since the admissible and inadmissible matrices can separately be
computed, we can parallelise the process to generate the 
${\mathcal H}$-matrix.
To achieve high parallel performance, however, we need to carefully
distribute the tasks to each thread because the size of each
block matrix may differ from each other.
In this study, the block matrices are stuck in descending order
in size and the block matrices are distributed to threads with dynamic
scheduling.
In the dynamic scheduling, all the threads share the list of tasks.
When we use $n$ threads for the parallelisation, $n$ tasks from top of
the list are firstly distributed to each thread.
When a thread finish its assigned task, the thread gets a new one from top of
the list of unprocessed tasks. By this scheduling, we can keep
almost all of the threads running and improve the parallel performance.

\subsection{Improvement of the ${\mathcal H}$-matrix structure by the agglomeration}
In the ${\mathcal H}$-matrix method, we sometimes are required to
refine a ${\mathcal H}$-matrix structure since the computational
cost for the ${\mathcal H}$LU will increase in accordance with 
fineness of the ${\mathcal H}$-matrix. The agglomeration is a technique
to refine a ${\mathcal H}$-matrix by approximating four matrix as one
low rank matrix. As an example, we consider the following matrix 
$A \in {\mathbb C}^{m\times n}$
constructed by four low rank matrices:
\begin{align}
 A=\left[
  \begin{array}{cc}
   U_1 V_1^{H} & U_3 V_3^{H} \\
   U_2 V_2^{H} & U_4 V_4^{H} 
  \end{array} 
 \right],
\end{align}
where $U_i\ (i=1,\cdots,4)$ (resp. $V_i$) denote a matrix whose
size is $m \times k_i$ (resp. $n \times k_i$) with the rank $k_i$.
In the agglomeration, we firstly express $A$ as a product of two
matrices $U \in {\mathbb C}^{m\times 2n}$ and 
$V \in {\mathbb C}^{n \times 2n}$ as follows:
\begin{align}
 A=U V^{H}
 =
 \left[
 \begin{array}{cccc}
  U_1 &     & U_3 & \\
      & U_2 &     & U_4 
 \end{array} 
 \right]
 \left[
 \begin{array}{cccc}
  V_1 & V_2 &     & \\
      &     & V_3 & V_4 
 \end{array} 
 \right]^{H}.
\end{align}
In order to approximate the matrix $A$ as a low rank matrix, we apply
the QR decomposition to $U$, $V$ and obtain the following expression
of $A$:
\begin{align}
 A=Q_U R_U R_V^{H} Q_V^{H}.
\end{align}
With the help of the singular value decomposition of $R_U R_V^{H}$,
$A$ is expressed in the form of $U^{'} \Sigma V^{'}$ in which 
$U^{'} \in {\mathbb C}^{m\times m}$, $V^{'} \in {\mathbb C}^{n\times n}$ 
denote orthogonal matrices, and $\Sigma \in {\mathbb C}^{m\times n}$
denote a diagonal matrix with singular values on the diagonal.
By rounding singular values less than a preset value, we approximate
$A$ as a low rank matrix.

\subsection{Acceleration of the sensitivity analysis with the ${\mathcal
 H}$-matrix method}
With the help of the ${\mathcal H}$-matrix method, we can reduce the
computational cost to solve the algebraic equations to $O(N\log N)$
for the number of boundary elements $N$ by the LU decomposition.
The accelerated LU decomposition with the ${\mathcal H}$-matrix method
is henceforth denoted as ${\mathcal H}$LU.
In this study, we further accelerate the sensitivity analysis by solving 
the forward problem (\ref{tmbvp1})--(\ref{tmbvp6}) 
and adjoint problem (\ref{adj1})--(\ref{adj2}) at one time with the
${\mathcal H}$LU.

Furthermore, we consider to accelerate the computation of
the electro-magnetic response
in the domain $\Omega_1$, $\Omega_2$ and adjoint incident field by
the fast matrix-vector multiplication with the ${\mathcal H}$-matrix method.
The TM polarised electro-magnetic response in each domain is summarised as follows:
\renewcommand{\arraystretch}{1.8}
\begin{alignat}{2}
 \left[
  \begin{array}{c}
   u(\bs x) \\
   \displaystyle
   \frac{\partial u(\bs x)}{\partial x_1} \\
   \displaystyle
   \frac{\partial u(\bs x)}{\partial x_2}
  \end{array}
 \right]
 &=
 \left[
  \begin{array}{ccc}
   \displaystyle
    \mu_1 S_1^{\Gamma_p}      & -D_1^{\Gamma_d}     & 
    \mu_1  S_1^{\Gamma_d} \\
   \mu_1 D_{1,1}^{\Gamma_p*} & -N_{1,1}^{\Gamma_d} & 
    \mu_1 D_{1,1}^{\Gamma_d*} \\
   \mu_1 D_{1,2}^{\Gamma_p*} & -N_{1,2}^{\Gamma_d} & 
    \mu_1 D_{1,2}^{\Gamma_d*}
  \end{array}
 \right]
 \left[
 \begin{array}{c}
  w_p \\
  u_d \\
  w_d 
 \end{array}
 \right]
 &&\ \ \ \bs x \in \Omega_1 ,
 \label{eq.inner1} \\
 \left[
  \begin{array}{c}
   u(\bs x) \\
   \displaystyle
   \frac{\partial u(\bs x)}{\partial x_1} \\
   \displaystyle
   \frac{\partial u(\bs x)}{\partial x_2}
  \end{array}
 \right]
 &=
 \left[
  \begin{array}{cc}
   \displaystyle
    D_2^{\Gamma_d}     & -\mu_2 S_2^{\Gamma_d} \\
    N_{2,1}^{\Gamma_d} & -\mu_2 D_{2,1}^{\Gamma_d*} \\
    N_{2,2}^{\Gamma_d} & -\mu_2 D_{2,2}^{\Gamma_d*}
  \end{array}
 \right]
 \left[
 \begin{array}{c}
  u_d \\
  w_d 
 \end{array}
 \right]
 &&\ \ \ \bs x \in \Omega_2 ,
 \label{eq.inner2}
\end{alignat}
where, $D_{i,j}^{\Gamma *}$ and $N_{i,j}^\Gamma$ respectively denote the following operators:
\begin{equation}
 [D_{i,j}^{\Gamma *} \psi](\bs x)=\int_\Gamma
  \frac{\partial G^i(\bs x,\bs y)}{\partial x_j}
  \psi(\bs y)
  ~\mathrm{d}\Gamma_y,
\end{equation}
\begin{equation}
 [N_{i,j}^\Gamma \phi](\bs x)=\int_\Gamma
  \frac{\partial^2 G^i(\bs x,\bs y)}{\partial x_j
  \partial n_y}
  \phi(\bs y)   ~\mathrm{d}\Gamma_y.
\end{equation}
The incident field of adjoint problem is also summarised as matrix
vector multiplication as follows:
\begin{align}
 \left[
  \begin{array}{c}
   \tilde {u}^\mathrm{inc}(\bs x_1) \\
   \vdots \\
   \tilde {u}^\mathrm{inc}(\bs x_n)
  \end{array}
 \right]
 =
 \left[
  \begin{array}{ccc}
   \displaystyle
    G^1(\bs x_1, \bs x_1^\mathrm{obs}) & \cdots & G^1(\bs x_1, \bs
    x_m^\mathrm{obs})\\
    \vdots & \ddots & \vdots \\
    G^1(\bs x_n, \bs x_1^\mathrm{obs}) & \cdots & G^1(\bs x_n, \bs
     x_m^\mathrm{obs}) 
  \end{array}
 \right]
 \left[
 \begin{array}{c}
  \displaystyle
  \frac{\partial f}{\partial u}(\bs x_1^\mathrm{obs}) \\
  \vdots \\
  \displaystyle
  \frac{\partial f}{\partial u}(\bs x_m^\mathrm{obs})  
 \end{array}
 \right]
 \ \ \ \bs x \in \Omega_1 ,
 \label{eq.adinc} 
\end{align}
in which $\bs x_i\ (i=1,\cdots,n)$ denotes the points at which we compute
the adjoint incident field.
We reduce the computational cost for assembling the coefficient matrices
and matrix vector multiplication with the ${\mathcal H}$-matrix method and
the ACA.
\section{Numerical examples}
In this section, we show the effectiveness of the proposed method
by some numerical examples.

\begin{figure}[htbp]
  \begin{center}
   \includegraphics[scale=0.5]{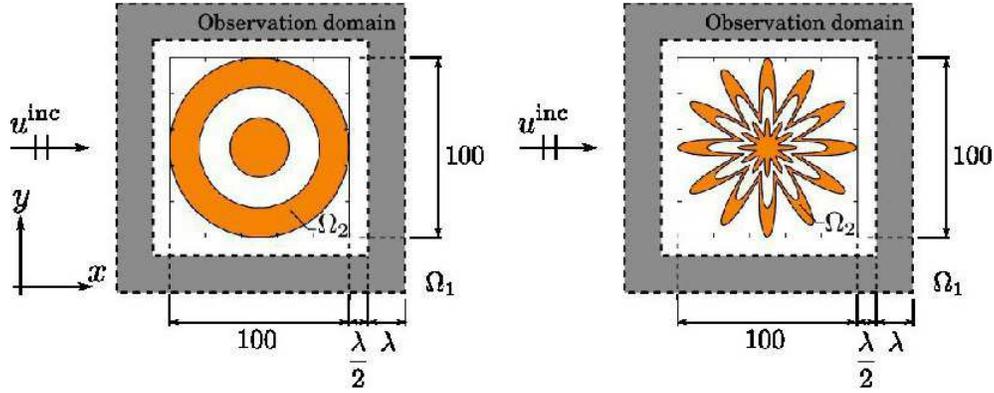}
   \caption{Problem statement for the test of performance of the 
   ${\mathcal H}$-matrix method.}
   \label{fig.prblm_une}
  \end{center}
\end{figure}

As the first numerical example, we investigate the computational cost
for the ${\mathcal H}$LU, computation of the electro-magnetic response in
the domain $\Omega_1 \cap D$, $\Omega_2$ and computation of the adjoint
incident field.
We consider a single step of an optimisation problem to minimise the amount of scattered
field at some observation points allocated in a lattice form around a
fixed design domain for an incident field which propagates in the
$x$-direction (\figref{fig.prblm_une}).
The number of observation points is $2290$, and the distance between
neighbouring observation points is fixed to $2.5$.
We compare the computational cost to solve the forward and adjoint
problems on the boundary of a dielectric elements shown in
\figref{fig.prblm_une} with the ${\mathcal H}$LU and compare it with that
for the GMRES accelerated by the FMM (FMGMRES) and GMRES without 
acceleration.
Note that since FMGMRES and GMRES are iterative solver, the forward and
adjoint problems are solved individually.
On the other hand, since the ${\mathcal H}$LU is a direct solver,
the LU decomposition of the coefficient matrix can be used for both
forward and adjoint analysis.

\begin{figure}[htbp]
  \begin{center}
   \includegraphics[scale=0.35]{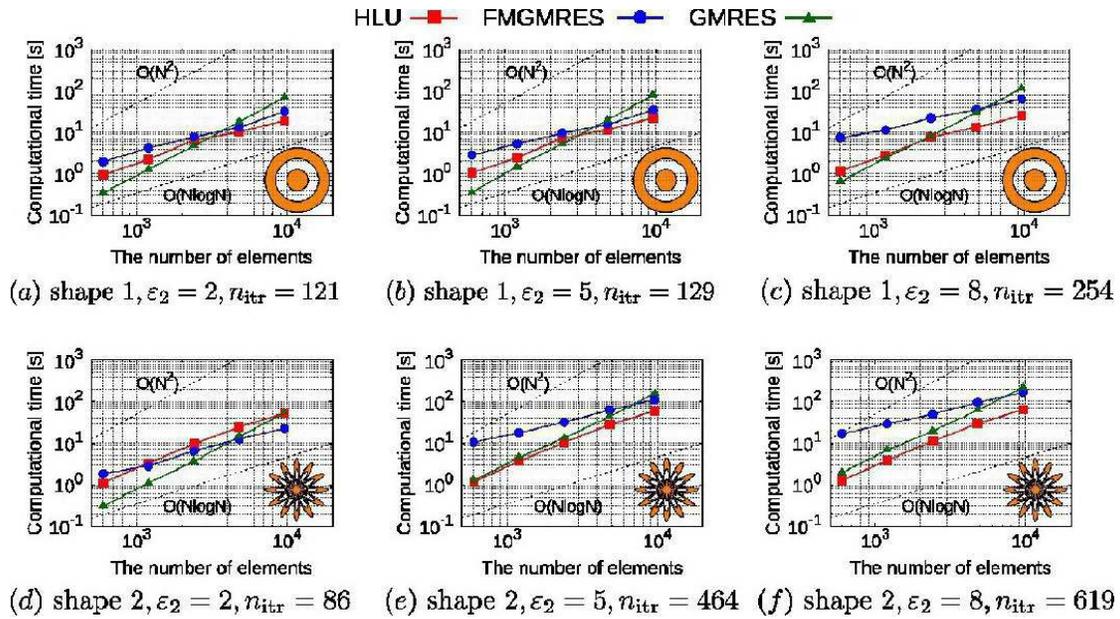}
   \caption{Computational time to compute the forward and adjoint
   response on the boundary of dielectric elements attached to each
   figure. Permittivity of dielectric elements $\varepsilon_2$ and the
   number of iteration $n_\mathrm{itr}$ for GMRES when $N$ equals to 
   $600$ are inserted.}
   \label{fig.time_une}
  \end{center}
\end{figure}

\begin{figure}[htbp]
  \begin{center}
   \includegraphics[scale=0.35]{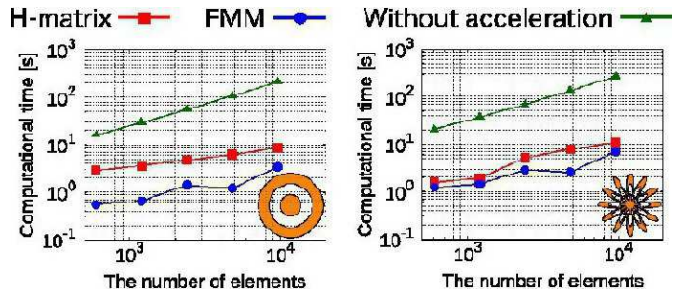}
   \caption{Computational time for the inner computation for
   $\varepsilon_2=2$.}
   \label{fig.time_inner}
  \end{center}
\end{figure}

The tolerance for ${\mathcal H}$-matrix 
operations and the ACA are set to $\varepsilon=10^{-5}$ in this study,
which is determined in a way such that relative error by the 
${\mathcal H}$-matrix method
for the solution on the boundary of a dielectric circle is negligible
compared to the discretisation error. Also, through some numerical
experiments, we employ $\eta=128$ and $n_\mathrm{min}=128$ which
reduce the computational cost for numerical analyses most efficiently.
Tolerance of the GMRES is set to $10^{-5}$. 
In every examples, we use the agglomeration technique when we generate a
${\mathcal H}$-matrix.
\figref{fig.time_une} shows the computational time for the number of
boundary elements $N = 600, 1200, 2400, 4800, 9600$ when $\varepsilon_2$
is either $2$ or $5$ or $8$. The number of iteration with the GMRES in
the case of $N=600$ is attached to each figure as $n_\mathrm{itr}$.
${\mathcal H}$LU achieves almost $O(N\log N)$ computational cost in every
cases.
Through the comparison between the result for different $\varepsilon_2$,
one observes that computational cost for FMGMRES and GMRES increase 
in association with increase of $\varepsilon_2$
due to deterioration of convergence property, 
while that for ${\mathcal H}$LU is almost the same independently on
$\varepsilon_2$.
As for the complexity of dielectric element shape, ${\mathcal H}$LU
shows stable computational time when $N$ is around $1000$.
On the other hand, for larger $N$ the computational time get to be
sensitive to the complexity of shape for all of solvers.
The computational time for the inner computation (\ref{eq.inner1})
and (\ref{eq.inner2}) in forward and adjoint problem for the same
parameters with the previous result is shown in \figref{fig.time_inner}.
Also, \tabref{tab.time_adjinc} shows the computational time for the
computation of the adjoint incident field on the boundary and inner
points for $N=9600$.
Both results indicate that the FMM is the fastest for a simple 
matrix-vector multiplication. 
The ${\mathcal H}$-matrix method also achieve the fast
computation compared to the case without acceleration.
One observes that the result for the FMM does not shows monotonous
increase in accordance with $N$.
This is because parameters for the FMM is chosen in a way that the
computational time for whole sensitivity analysis is the fastest.

\begin{table}[h!] 
 \begin{center}
  \caption{Computational time for the adjoint incident field for $\varepsilon_2=2$.}
  \label{tab.time_adjinc}
   \begin{tabular}{|c|c|c|}\hline
                             & Shape 1  &  Shape 2  \\ \hline
    ${\mathcal H}$-matrix method & 0.45 [s] &  0.40 [s] \\ \hline
    FMM                      & 0.13 [s] &  0.18 [s] \\ \hline
    Without acceleration     & 9.80 [s] &  9.44 [s] \\ \hline
   \end{tabular}
 \end{center} 
\end{table}

\begin{figure}[htbp]
  \begin{center}
    \includegraphics[scale=0.45]{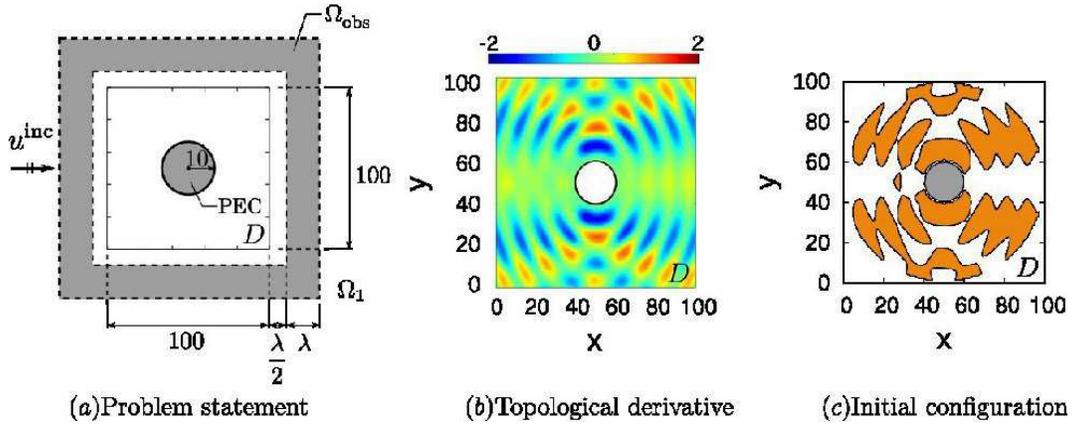}
   \caption{(a) Conventional problem statement, (b) The topological
   derivative when only circular PEC is allocated in the design
   domain, (c) The initial configuration determined by the sign of the
   topological derivative (b) at each point.}
    \label{fig.prblm_pec_tpo}
  \end{center}
\end{figure}

Next, we apply the proposed method to the conventional topology
optimisation problem (\ref{eq.conv_obj})--(\ref{tmbvp6}) of cloaking
devices.
We consider to determine configuration of cloaking
devices which make a circular PEC
invisible for a TM polarised plane incident field.
The radius of PEC is fixed to $10.0$ in this example.
We consider a design domain $D$ whose size is $[0,100]\otimes [0,100]$ 
(\figref{fig.prblm_pec_tpo}).
We allocate $2290$ observation points 
$\bs x_m^\mathrm{obs}\ (m=1,\cdots,2290)$ on each lattice points 
in an observation domain $\Omega_\mathrm{obs}$ around $D$ and define the
objective
function by equation (\ref{eq.conv_obj}) as sum of the scattered field 
on the observation points.
We firstly determine the initial configuration $\Omega_2^\mathrm{init}$
as follows:
\begin{align}
 \Omega_2^\mathrm{init}
 =
 \left\{
 \bs x ~|~ 
 {\mathcal T}^0(\bs x) \le 0,
 \left\|\bs x-(50, 50)^\mathrm{T} \right\|_2
 \le 50
 \right\}, \label{init_conf}
\end{align}
where ${\mathcal T}^0(\bs x)$ denotes the topological derivative when
only circular PEC is allocated in the design domain
(\ref{fig.prblm_pec_tpo}).

\begin{figure}[htbp]
  \begin{center}
   \includegraphics[scale=0.8]{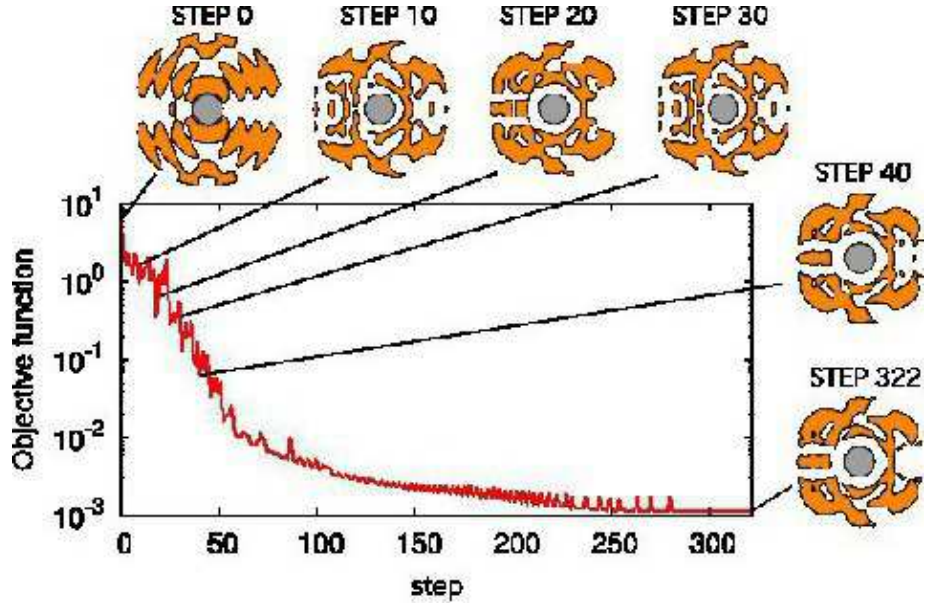}
    \caption{The history of the objective function and configuration of
   dielectric elements for $\varepsilon_2 = 2$, 
   $\tau = 5.0\times 10^{-3}$.} \label{fig.obj_hist}
  \end{center}
\end{figure}

\begin{figure}[htbp]
  \begin{center}
   \includegraphics[scale=0.7]{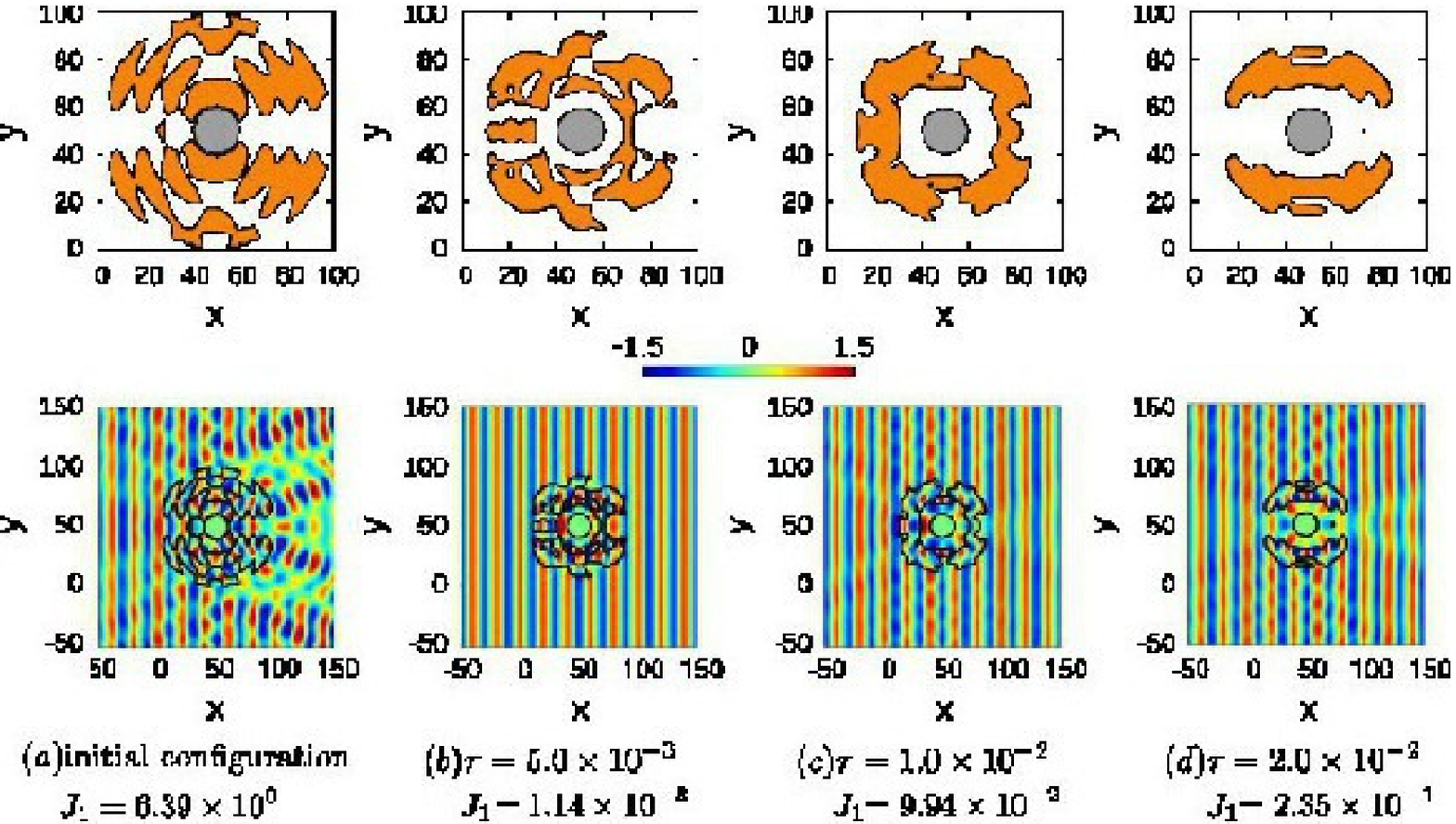}
    \caption{The initial configuration (a) and the optimal configuration
   (b)--(c) with the electric responses for the conventional
   optimisation problem in the case of $\varepsilon_2=2$.}
   \label{fig.cloak_pec1_eps2}
  \end{center}
\end{figure}

\begin{figure}[htbp]
  \begin{center}
    \includegraphics[scale=0.7]{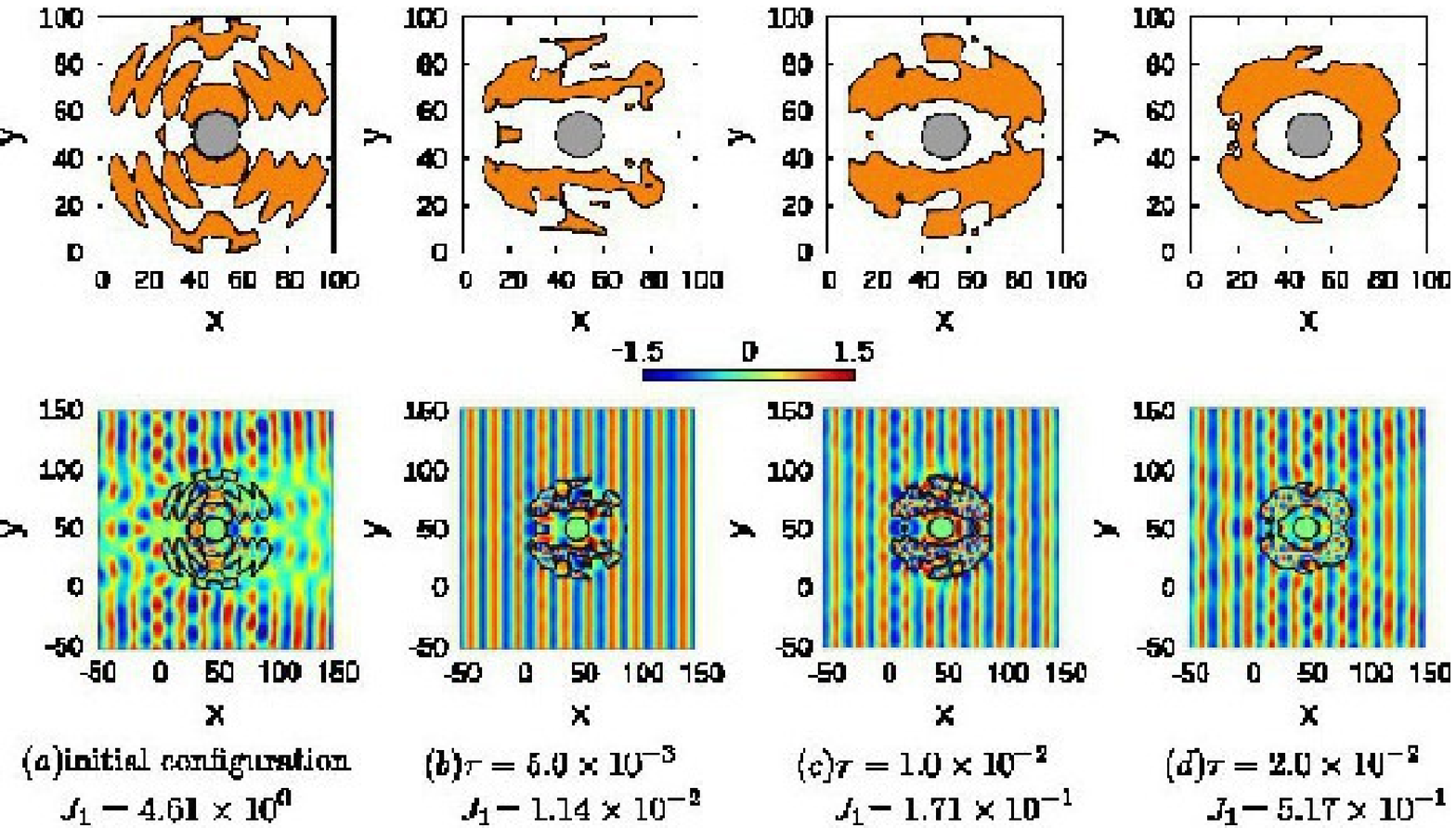}
    \caption{The initial configuration (a) and the optimal configuration
   (b)--(c) with the electric responses for the conventional
   optimisation problem in the case of $\varepsilon_2=5$.}
   \label{fig.cloak_pec1_eps5}   
  \end{center}
\end{figure}

\begin{figure}[htbp]
  \begin{center}
    \includegraphics[scale=0.55]{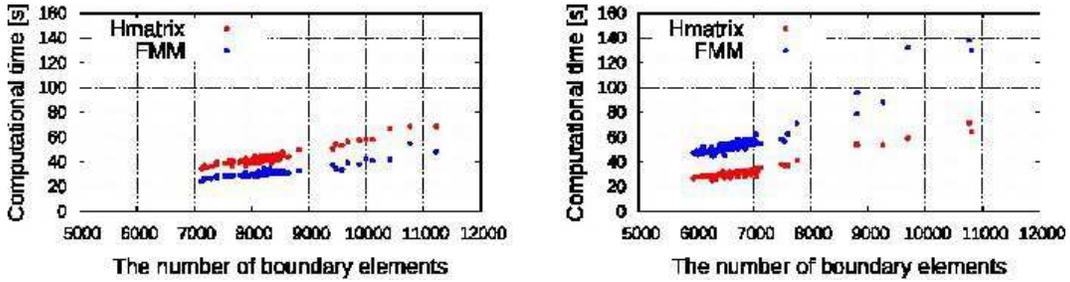}
   \caption{Computational time for the sensitivity analysis at each
   step of the conventional optimisation problem for $\tau=5.0\times
   10^{-3}$, $\varepsilon_2=2$ (left) and $5$ (right).}
    \label{fig.time_cloak_pec1}
  \end{center}
\end{figure}

In the following results, the value of the objective function is
normalised by the one when only circular PEC is allocated, which
is denoted as $J_1$.
The history of the objective function for $\tau = 5.0\times 10^{-3}$ 
and $\varepsilon_2=2.0$ is shown in \figref{fig.obj_hist}.
The objective function is successfully reduced as the optimisation step
goes.
For the optimal configuration, 
the objective function is reduced to approximately $0.11\%$ of the
original objective function when only circular PEC is allocated.
The initial shape, optimal configurations and the
electric field for 
$\tau=5.0\times 10^{-3}, 1.0\times 10^{-2}, 2.0\times 10^{-2}$
, $\varepsilon_2=2$ (resp. $\varepsilon_2=5$) are shown in
\figref{fig.cloak_pec1_eps2} (resp. \figref{fig.cloak_pec1_eps5}).
For the initial configuration, the electric field is highly affected
by scattered field. On the other hand, the optimal configurations for 
$\tau=5.0\times 10^{-3}$ successfully reduce the scattering around the
PEC and dielectric elements.
One observes that complexity of the optimal configuration is reduced
corresponding to the increase of $\tau$ while the objective function
for the optimal configuration becomes larger.
\figref{fig.time_cloak_pec1} shows comparison of the computational time
for the sensitivity analysis with the ${\mathcal H}$LU
and that with the FMGMRES at each optimisation step.
One observes that computational time for FMGMRES is sensitive to
the change of the permittivity $\varepsilon_2$, while the ${\mathcal H}$LU
takes almost the same cost for the sensitivity analysis 
independently on $\varepsilon_2$.
Total computational time to obtain the optimal configuration with the
${\mathcal H}$LU for $\varepsilon_2=5$ is reduced to $58.7\%$ of that with 
the FMGMRES.

\begin{figure}[htbp]
  \begin{center}
    \includegraphics[scale=0.5]{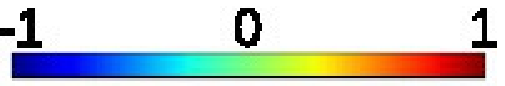}
    \includegraphics[scale=0.7]{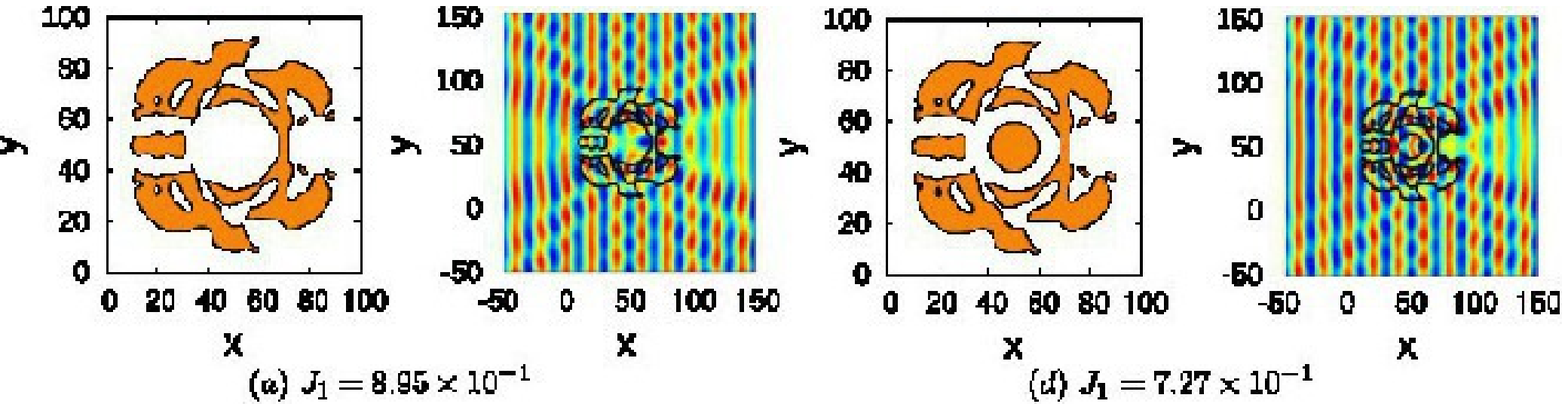}
    \includegraphics[scale=0.7]{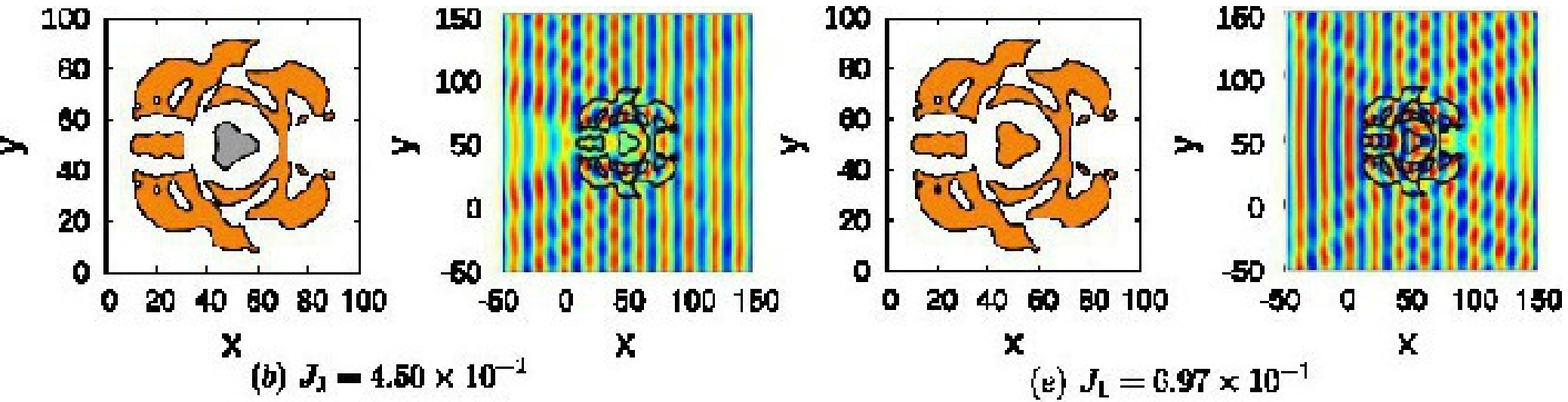}
   \includegraphics[scale=0.7]{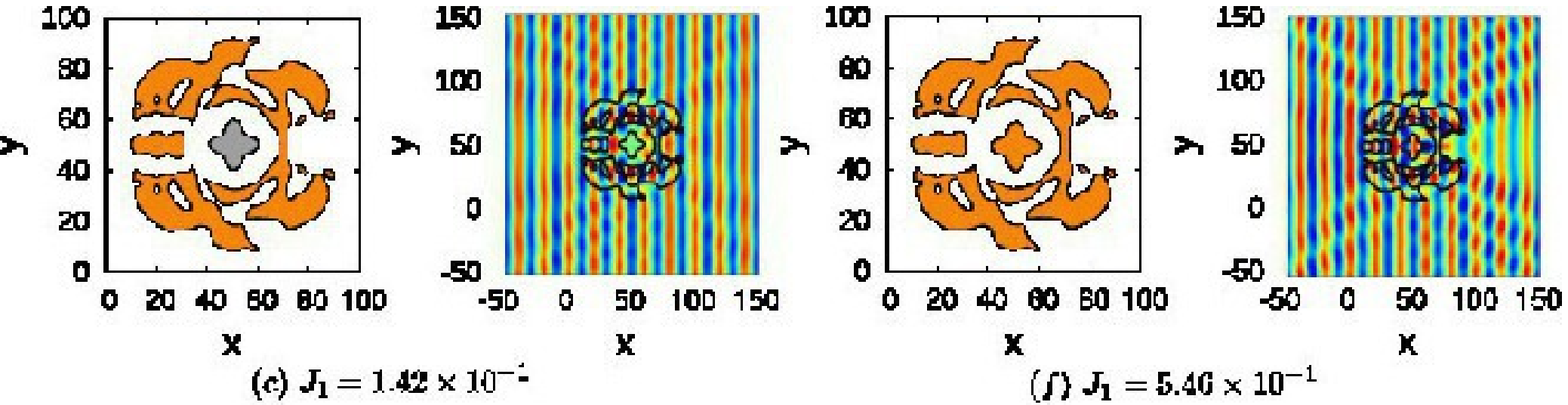}
    \caption{The electric field when nothing is put (a), some shapes
   of PEC are put (b), (c) and some shapes of dielectric element are put
   (d)--(f) in the cloaking device obtained by
   the conventional optimisation problem with $\varepsilon_2=2$,
   $\tau=5.0\times 10^{-3}$.}
   \label{fig.cloak_pec1_pecshape}
  \end{center}
\end{figure}

In order to check the dependency of the cloaking effect on shape and
material of a target object,
we show the distribution of the electric field when various shape of PEC
or dielectric material is allocated in the obtained cloaking device
for $\varepsilon_2=2$, $\tau=5.0\times 10^{-3}$
(\figref{fig.cloak_pec1_pecshape}).
One observes that the value of the
objective function differs by the shape and material of the target objects,
and obtained design does not work as a cloaking device for target
objects except for circular PEC.

\begin{figure}[htbp]
  \begin{center}
    \includegraphics[scale=0.45]{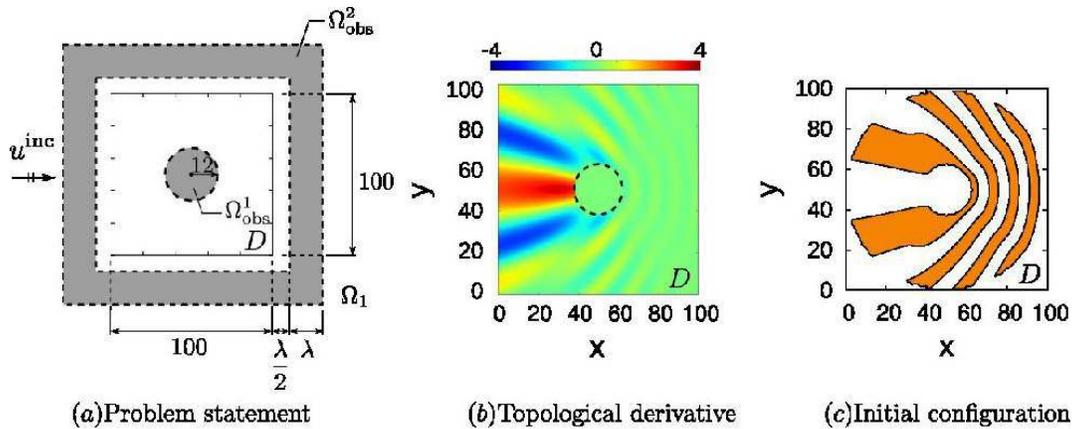}
   \caption{(a) Modified problem statement, (b) The topological derivative when
   nothing is allocated in the design domain, (c) The initial
   configuration determined by the sign of the topological
   derivative (b) at each point.}
    \label{fig.prblm_tpo}
  \end{center}
\end{figure}

In the next example, we redefine the objective function by
(\ref{mod_obj}) as sum of the electric field in a fixed circle
$\Omega_\mathrm{obs}^1$ whose
radius is $12.0$ at the center of the design domain and the scattered
field in $\Omega_\mathrm{obs}^2$ which is the same domain with
$\Omega_\mathrm{obs}$ in the previous example (\figref{fig.prblm_tpo}).
We put $437$ observation points on lattice points in
$\Omega_\mathrm{obs}^1$ in a way that
the distance between neighbouring points is $1$ in addition
to the $2290$ observation points in $\Omega_\mathrm{obs}^2$.
The initial configuration is obtained by redefining ${\mathcal T}^0(\bs x)$
in the equation (\ref{init_conf})
as the topological derivative when nothing is put in the design domain
(\figref{fig.prblm_tpo}).

\begin{figure}[htbp]
  \begin{center}
   \includegraphics[scale=0.7]{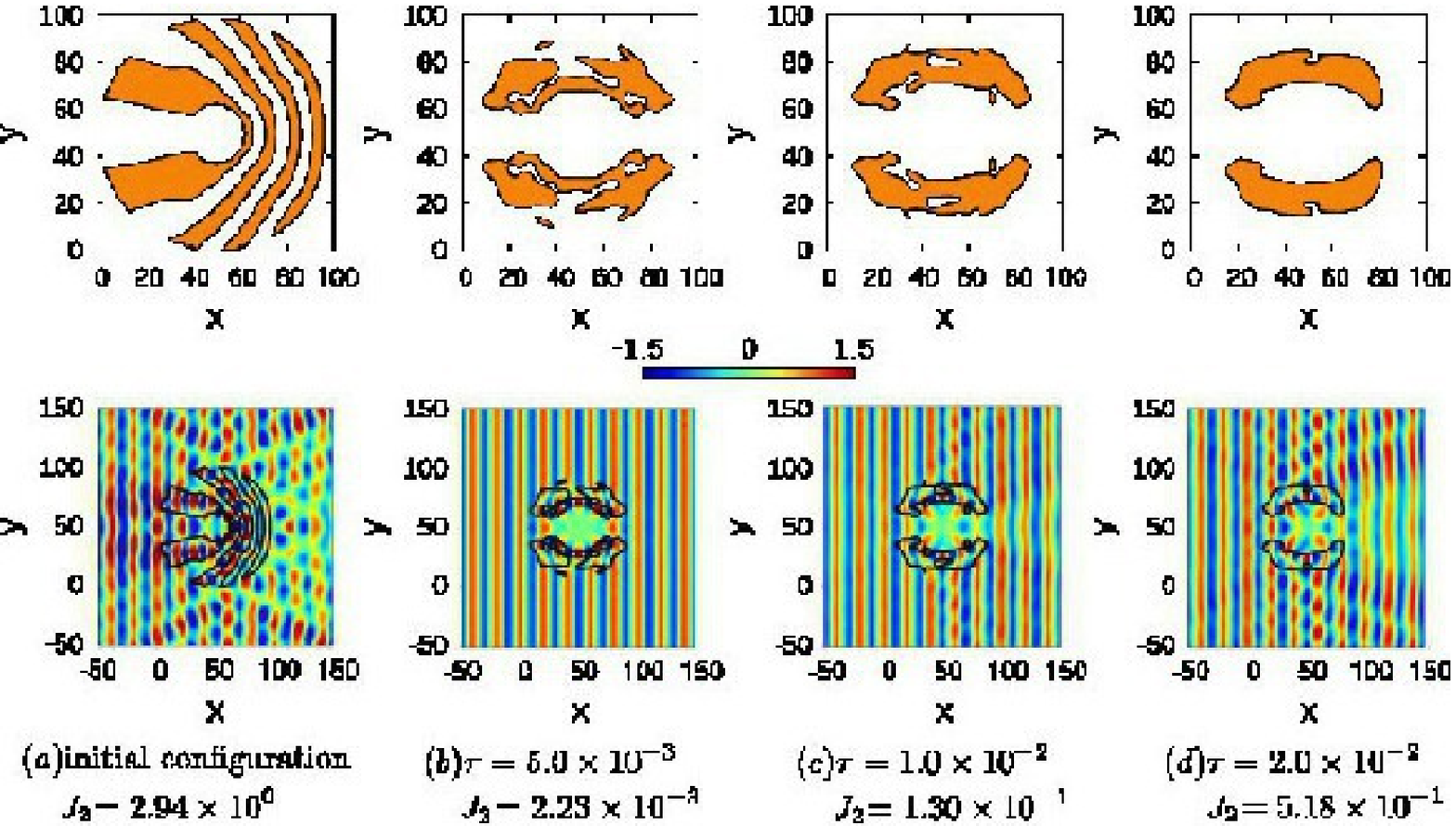}
   \caption{The initial configuration (a) and the optimal configuration
   (b)--(c) with the electric responses for the modified
   optimisation problem in the case of $\varepsilon_2=2$.}
    \label{fig.cloak_pec0_eps2}
  \end{center}
\end{figure}

\begin{figure}[htbp]
  \begin{center}
   \includegraphics[scale=0.7]{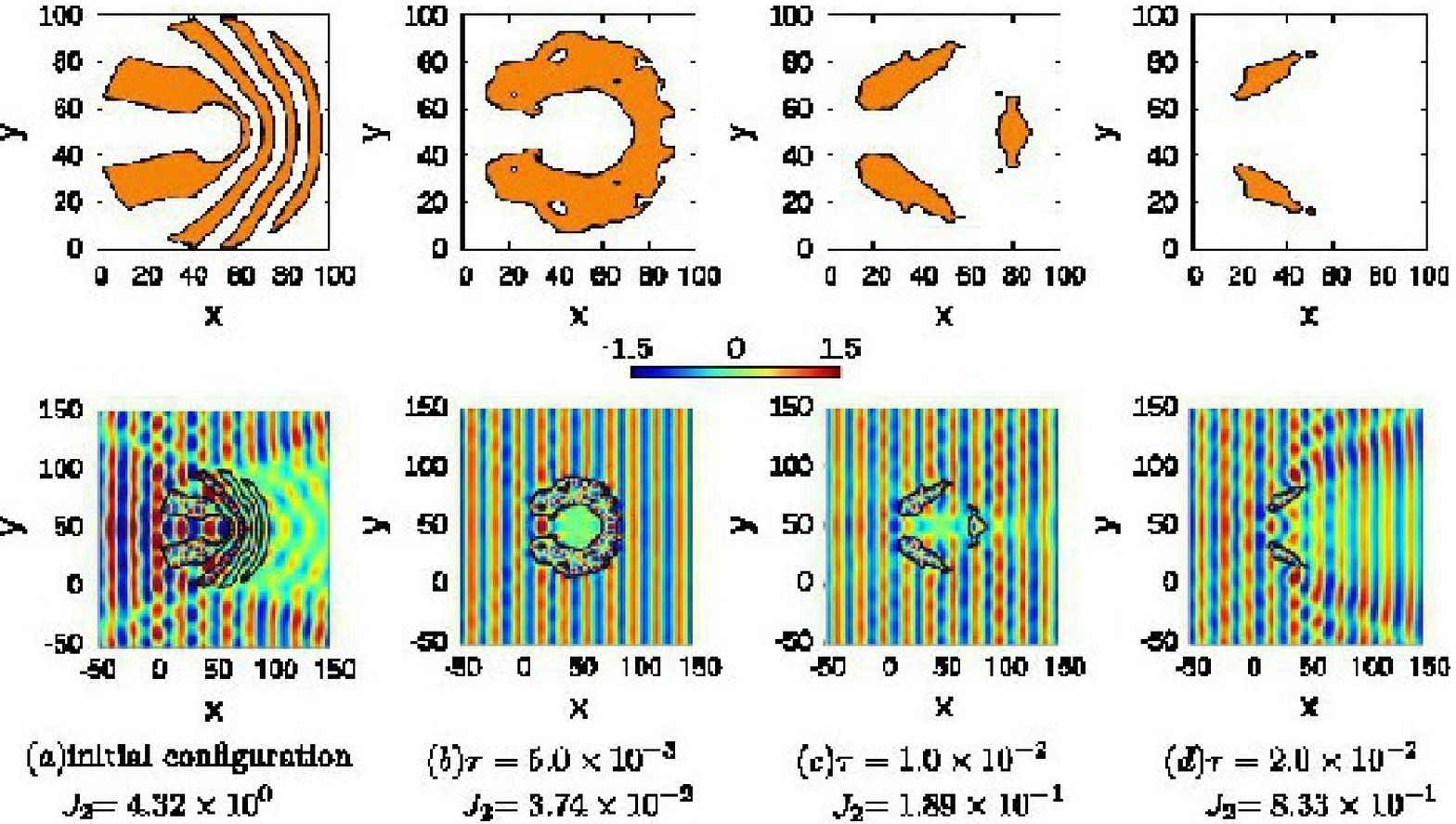}
   \caption{The initial configuration (a) and the optimal configuration
   (b)--(c) with the electric responses for the modified
   optimisation problem in the case of $\varepsilon_2=5$.}
    \label{fig.cloak_pec0_eps5}
  \end{center}
\end{figure}

\begin{figure}[htbp]
  \begin{center}
   \includegraphics[scale=0.55]{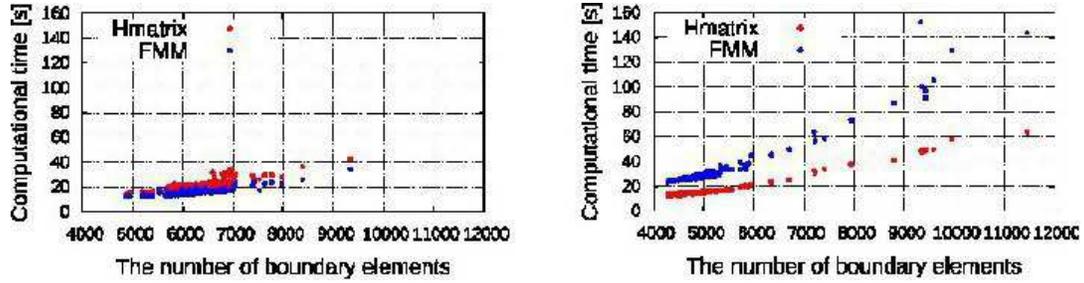}
   \caption{Computational time for the sensitivity analysis at each
   step of the modified optimisation problem for $\tau=5.0\times
   10^{-3}$, $\varepsilon_2=2$ (left) and $5$ (right).}
   \label{fig.time_cloak_pec0}
  \end{center}
\end{figure}

\begin{figure}[htbp]
  \begin{center}
    \includegraphics[scale=0.5]{fig/colorbar_1.eps}
    \includegraphics[scale=0.7]{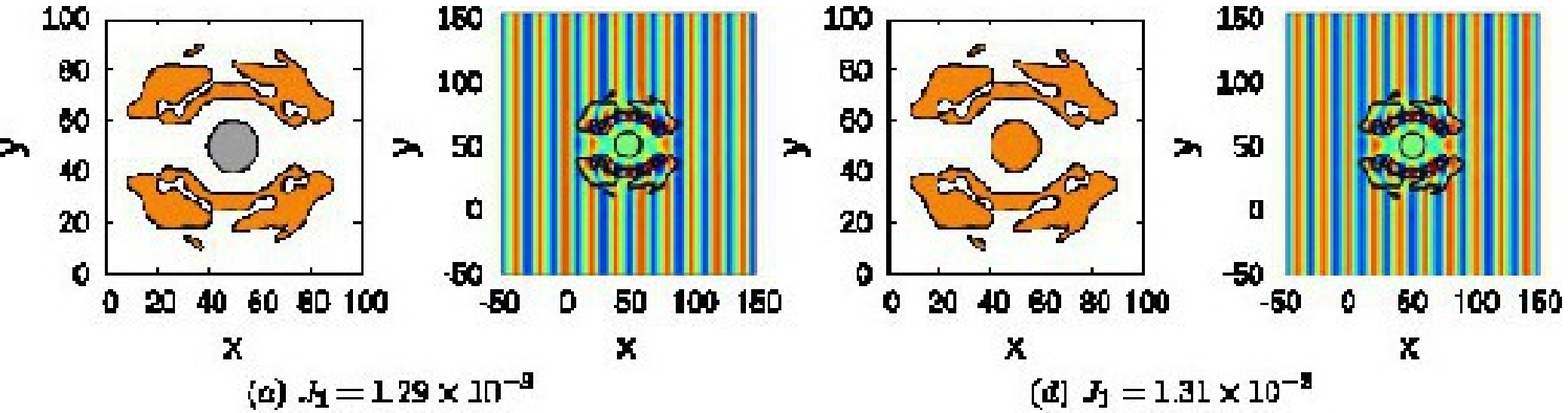}
    \includegraphics[scale=0.7]{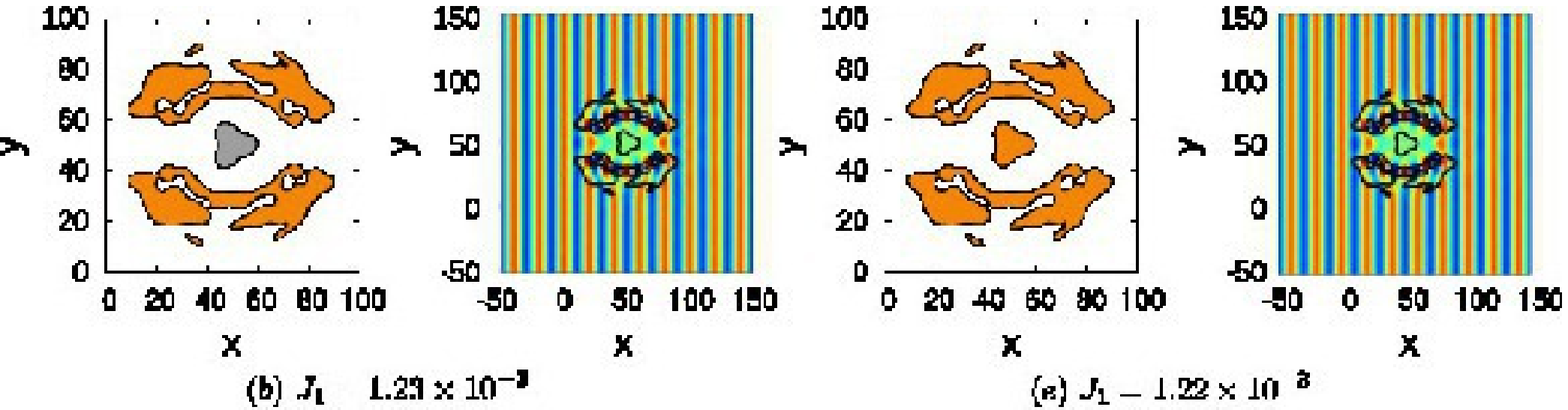}
   \includegraphics[scale=0.7]{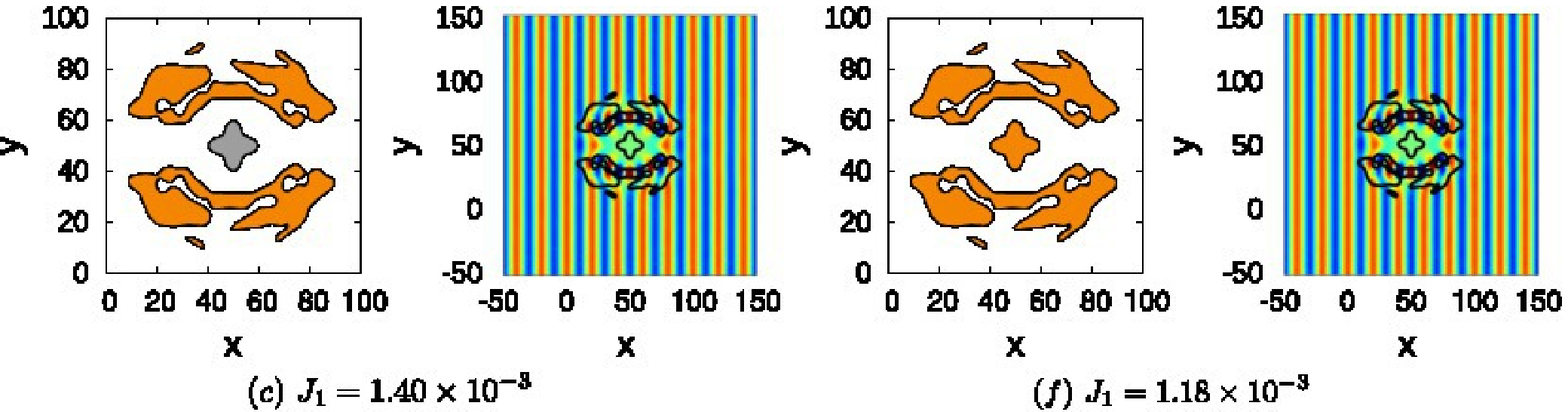}
    \caption{The electric field when nothing is put (a), some shapes
   of PEC are put (b), (c) and some shapes of dielectric element are put
   (d)--(f) in the cloaking device obtained by
   the modified optimisation problem with $\varepsilon_2=2$,
   $\tau=5.0\times 10^{-3}$.}
   \label{fig.cloak_pec0_pecshape}
  \end{center}
\end{figure}

In the following results, we denote the objective function which is
normalised by the one when nothing is put in the design domain as $J_2$.
Obtained configuration of dielectric elements and distribution of the
electric field for $\varepsilon_2=2$ and $5$ are shown in
\figref{fig.cloak_pec0_eps2} and \figref{fig.cloak_pec0_eps5},
respectively.
In both cases,
the optimal configuration for $\tau=5.0\times 10^{-3}$ successfully
reduces the electric field in $\Omega_\mathrm{obs}^1$ in addition to
the scatterred field in $\Omega_\mathrm{obs}^2$.
As with the previous example, one finds that higher value of $\tau$ leads
more simple configuration and larger objective function.
Computational time for the sensitivity analysis at each optimisation
step (\figref{fig.time_cloak_pec0}) shows that the ${\mathcal H}$LU works
stably independently on $\varepsilon_2$ compared to the FMGMRES.
\figref{fig.cloak_pec0_pecshape} shows the electric response of the optimal
configuration for $\varepsilon_2=5$ and $\tau=5.0\times 10^{-3}$
when various shape of PEC or dielectric element is allocated in
$\Omega_\mathrm{obs}^2$.
For the comparison with the previous result
\figref{fig.cloak_pec1_pecshape}, the value of the objective function $J_1$
defined by the conventional manner is attached.
The result indicates that the obtained cloaking design works
successfully independently on the shape and material of the hidden object.



\section{Conclusion}
We developed the topology optimisation method of cloaking devices
which work for arbitrary-shaped target object with efficient and accurate
sensitivity analysis using the BEM and the ${\mathcal H}$-matrix method. 
By the proposed method, we obtained structure of
cloaking devices with low computational cost and confirmed that it works
successfully independently on the shape of perfect electric conductor
allocated in the cloaking device. As a future task, we consider to
extend the proposed method to the 3D problem. In 3D cases, the numerical
cost for the electro-magnetic field analysis will be a key factor for
the topology optimisation. Hence, we firstly try to develop a fast 
BEM with the ${\mathcal H}$-matrix method for the Maxwell's equations.
Also, in 3D problems, the memory consumption to memorise the matrix
for the inner computation will be huge.
Hence, we need to develop an efficient algorithm to reduce the memory
consumption.

\bibliographystyle{abbrv}
\bibliography{ref} 
\end{document}